\documentclass[english,letter,11pt,twosided]{article}

\usepackage{amsfonts}
\usepackage[left=2.3cm,right=2.0cm, bottom = 2.3cm, top=2.4cm]{geometry}
\usepackage[utf8]{inputenc}
\usepackage{babel}
\usepackage{slantsc}
\usepackage{array}
\usepackage{amsmath}
\usepackage{amsthm, mathtools}
\usepackage{amssymb}
\usepackage[pdftex]{color,graphicx}
\usepackage{enumerate}
\usepackage{mathtools}
\usepackage{bbm}
\usepackage{tikz}
\usepackage{subfigure}
\usepackage{pgfplots}
\usepackage{utopia}
\usepackage[labelsep=period]{caption}
\usepackage{hyperref}
\usepackage{setspace}

\usetikzlibrary{patterns,snakes}
\newtheorem{theorem}{Theorem}

\newtheorem{corollary}{Corollary}

\theoremstyle{definition}
\newtheorem{definition}{Definition}

\newtheorem{assumption}{Assumption}
\newtheorem{remark}{Remark}

\newcommand{\ignore}[1]{}

\newcommand{\R}{{\mathbb{R}}}
\newcommand{\Q}{{\mathcal{Q}}}

\newcommand{\G}{{\mathcal{G}}}
\newcommand{\K}{{\mathcal{K}}}
\newcommand{\V}{{\mathcal{V}}}
\newcommand{\E}{{\mathcal{E}}}
\newcommand{\M}{{\mathcal{M}}}

\newcommand{\ie}{{i.e., }}

\newcommand{\Tr}{\mathbf{Tr}}
\newcommand{\adj}{A_{\G}}
\newcommand{\hs}{\hspace{0.05cm}}

\newcommand{\NN}{\boldsymbol{\mathcal N}}
\newcommand{\RR}{\boldsymbol{\rho}_{\textbf{ss}}}

\definecolor{PennBlue}{RGB}{001,031,091}
\definecolor{PennRed}{RGB}{153,0,0}
\definecolor{NewBlue}{RGB}{001,031,110}
\definecolor{NewRed}{RGB}{200,0,0}
\hypersetup{
pdfborder = {0 0 0},
    colorlinks,
    citecolor=PennRed,
    filecolor=PennRed,
    linkcolor=PennRed,
    urlcolor=PennRed
}
\begin{document}

\title{\fontsize{21}{21}%
\selectfont {Fundamental Limits and Tradeoffs on  Disturbance Propagation in Large-Scale  Dynamical Networks}\thanks{{\footnotesize This work was supported by the Office of Naval Research under award ONR N00014-13-1-0636.}}}
\author{\fontsize{13}{13} {Milad Siami}\thanks{{\footnotesize Lehigh University, Bethlehem, PA, USA.}}  \hspace{0.5in}{Nader Motee}\thanks{{\footnotesize Lehigh University, Bethlehem, PA, USA.}}  \\
~}
\date{\fontsize{11}{11}\selectfont This version: December 2014 \\
First version: March 2014} 
\maketitle

\begin{abstract}
\fontsize{10}{10}\selectfont \baselineskip0.5cm
We consider performance deterioration  of interconnected linear dynamical networks subject to exogenous  stochastic disturbances. The focus of this paper is on first-order and second-order linear consensus networks. We employ the expected value of the steady state dispersion of the state of the entire network as a performance measure and develop a graph-theoretic methodology to relate structural specifications of the underlying graphs  of the network to the performance measure. We explicitly quantify several inherent fundamental limits on the best achievable levels of performance in linear consensus networks and show that these limits of performance are merely imposed by the specific structure of the underlying graphs. Furthermore, we discover new connections between notions of sparsity and the performance measure. Particularly, we characterize several fundamental tradeoffs that reveal interplay  between the performance measure and various sparsity measures of a linear consensus network.  At the end, we apply our results to two real-world dynamical networks and provide energy interpretations for the proposed performance measures. It is shown that the total power loss in synchronous power networks and total kinetic energy of a network of autonomous vehicles in a formation are viable performance measure for these networks and fundamental limits on these  measures quantify the best achievable levels of energy-efficiency in these dynamical networks. 
\textit{IEEE Transaction on Automatic Control, Under Review.}  \newline
\newline
\end{abstract}
\fontsize{12}{12}\selectfont


\thispagestyle{empty} \newpage \fontsize{11}{11}\selectfont\baselineskip %
0.60cm

\onehalfspacing

\section{Introduction}
	The issue of fundamental limits and their corresponding tradeoffs in large-scale interconnected dynamical systems design lie at the very core of theory of distributed feedback control systems as it reveals what is achievable and conversely what is not achievable by distributed feedback control laws. Improving global performance as well as robustness to exogenous disturbances in dynamical networks are crucial for sustainability in engineered infrastructures; examples include a group of autonomous vehicles such as UAVs in a formation, distributed emergency response systems, interconnected transportation networks, energy and power networks, metabolic pathways and even social networks \cite{Bamieh12,Barooah,Siami13siam,Siami13acc, Siami13cdc, MoteeCBKD10,abbas,Barooah-acc-12, Scardovi2010, Zelazo2013}. One of the outstanding analysis problems in dynamical networks is to devise a mathematical methodology to study and characterize intrinsic fundamental limits and their tradeoffs in networks of interconnected systems. Providing solutions to  this important challenge will enable us to develop underpinning principles to design robust-by-design dynamical networks that are less fragile to exogenous disturbances.

	There have been several recent works on the performance analysis of first-order and second-order linear consensus networks; only to name a few, we refer to \cite{Bamieh12, Siami13cdc, abbas, Young10,Bamieh11,  Zelazo-Mesbahi, LovisariGarinZampieriResistance, Nicola, Lin2014, Jadbabaie13} and references in there.  The reference papers \cite{Bamieh12, Young10, Bamieh11, Jadbabaie13} study performance of a class of linear consensus networks under influence of white exogenous noises. The common approach of the  above-mentioned papers is to adopt the $\mathcal H_2$-norm of the system (from the disturbance input to the performance output of the system) as a scalar performance measure. Since Laplacian matrices belong to the class of normal matrices, the $\mathcal H_2$-norm of linear consensus networks can be exactly calculated  as a function of the eigenvalues of the state matrix of the system \cite{Bamieh12}. When the state matrix of the system is a graph Laplacian, this scalar measure is proportional to the total effective resistance of the graph. The concept of effective resistance has been used in several disciplines and applications. In the context of electric circuit analysis, the effective resistance of an edge is the resistance measured between endpoints of that edge. In the context of random walks and Markov chains on networks, the effective resistance of an edge can be interpreted as the commute time between the endpoints of that edge.  Another interesting version of the notion of effective resistance appears in the context of graph sparsification, where the goal is to approximate a given graph by a sparse graph. In this setting, the effective resistance is defined as probability of appearing an edge in a random spanning tree of the graph (see \cite{Spielman} and references in there). In \cite{Barooah2006}, the authors demonstrate a physical interpretation of the effective resistance in least-squares estimations  as well as motion control problems.

	Besides the effective resistance interpretations of $\mathcal H_2$-norm of a linear consensus network, this measure can be viewed as a macroscopic performance measure that captures the notion of coherence in large-scale dynamical networks. In \cite{Bamieh12}, linear consensus networks over multi-dimensional discrete torus coupling graphs are considered and it is shown that how the $\mathcal H_2$-norm of such networks scale asymptotically  with the network size. In \cite{Zelazo-Mesbahi}, the authors consider the $\mathcal H_2$-norm performance measure for a class of first-order consensus networks with exogenous inputs in the form of process and sensor noises. The performance measure used in \cite{Zelazo-Mesbahi} is different from those scalar measures considered in \cite{Bamieh12, Young10, Siami13cdc, Siami13siam}. The proposed analysis method in \cite{Zelazo-Mesbahi} employes the edge agreement protocol by considering a minimal realization of the edge interpretation system. Another related work is reported in \cite{Jadbabaie13}, where the authors use the Euclidean norm coefficient of ergodicity to find upper bounds on the $\mathcal H_2$-norm performance measure.

In this paper, we study first-order and second-order linear consensus algorithms for large-scale dynamical networks that are subject to exogenous stochastic disturbance inputs. In Section \ref{sec:first-order}, the steady state variance of the output of the dynamical network is employed as a performance measure in order to quantify to what degrees the performance of the network deteriorates as the result of disturbance input. It can be shown that this performance measure is equal to the square of the $\mathcal{H}_2$-norm of the network from the disturbance input to the output. This performance measure has an interesting output energy interpretation if we consider a linear consensus network with identically zero  input and perturb the network by a random initial condition. The value of the proposed performance measures is equal to the average output energy needed to be consumed throughout the network in order to steer the state of  this randomly perturbed dynamical network  to its equilibrium point. The physical interpretation of this performance measure depends on the application and is domain-specific. We show in Section \ref{sec:applications} that the total resistive power loss in a linearized model of interconnected network of synchronous generators and the total kinetic energy (also known as flock energy) of a group of controlled vehicles in a formation are examples of  admissible performance measures according to our definition. 

	Our primary focus  is to highlight the important role of underlying graphs of linear dynamical networks in emergence of severe theoretical hard limits on the best achievable levels of global performance. The structure of the underlying {\it coupling graph} of a dynamical network depends on the coupling structure among the subsystems, which are usually imposed by physical laws and/or global objectives. We consider linear time-invariant networks that are operating in closed-loop, i.e., linear dynamical networks that have been already  stabilized by a linear state feedback control law. In some applications such as formation control of autonomous vehicles, sparsity pattern of the underlying information structure in the controller array determines communication requirements among the subsystems, and  as a result, it defines the sparsity pattern of the underlying coupling graph of the closed-loop network. 
	
The first contribution of this paper is quantification of inherent fundamental limits on the best achievable values for the performance measure for linear consensus  networks. Our results in Section \ref{sec-IV} are classified with respect to unweighted and weighted underlying coupling graphs. It is shown that the performance measure of first-order consensus networks is\footnote{We employ  the big omega notation in order to generalize the concept of asymptotic lower bound in the same
way as $\mathcal{O}$ generalizes the concept of asymptotic upper bound. We adopt the following definition according to \cite{Knuth76}:
\begin{equation} 
f(n) = \Omega(g(n)) ~\Leftrightarrow~g(n) = \mathcal{O}(f(n)), \label{big-omega}
\end{equation}
where $\mathcal{O}$ represents the big O notation. In the left hand side of \eqref{big-omega}, the $\Omega$ notation implies that $f(n)$ grows at least of the order of $g(n)$.} $\Omega(n)$ for networks with fairly sparse  unweighted coupling graphs such as tree and unicyclic graphs, where $n$ is the network size. It is $\Omega(1)$ for networks with fairly dense unweighted coupling graphs such as complete bipartite and complete graphs. The performance measure is $\mathcal{O}(n^2)$, where networks with path-like coupling graphs experience the worst levels of performance.  The performance measure of a (Type 2) second-order consensus network is $\Omega(n^{-1})$ for networks with fairly dense and $\Omega(n)$ for networks with fairly sparse  unweighted coupling graphs; furthermore, it is $\mathcal{O}(n^4)$. For linear  consensus networks with weighted coupling graphs, it is shown that by subsuming more graph specifications in our calculations one can obtain improved lower bounds for the best achievable values for the performance measure. Our extensive simulation results assert that our theoretical lower bounds are tighter for networks with  rather dense coupling graphs. 

The impacts of presented fundamental limits in Section \ref{sec-IV} usually appear as critical interplays between various performance   and sparsity measures in consensus networks. Our second contribution is characterization of  several intrinsic tradeoffs between sparsity features of a coupling graph and global performance measures of linear consensus networks. In Section \ref{sec-V}, we formulate several uncertainty-principle-like inequalities that assert that networks with more sparse coupling graphs incur higher values of the performance measure. 

In Section \ref{sec:applications}, we consider two real-world networks and compute their corresponding performance measures.  We show that the total power loss in synchronous power networks and total kinetic energy of a network of autonomous vehicles in a formation are viable performance measure for these networks. The interpretation of our theoretical results in Section \ref{sec-IV} indicates that  existence of  fundamental limits on these performance measures quantify the best achievable levels of energy-efficiency in these dynamical networks.

\begin{remark}
The proof of all theorems and corollaries are given in the appendix at the end of the paper. 
\end{remark}


\section{Mathematical Notations}

{\it Matrix Theory:}	The set of all nonnegative real numbers is denoted by $\mathbb R _+$. The $n \times 1$  vector of all ones is denoted by $\mathbf{1}_{n}$, the $n \times n$ identity matrix by $I_{n}$, the $m \times n$ zero matrix by $\mathbf 0_{m \times n}$, and the $n \times n$ matrix of all ones by $J_{n}$. We will eliminate subindices of these matrices whenever the appropriate dimensions are clear from the context. The {\it centering matrix} of size $n$ is defined by
\begin{equation*}
M_n := I_{n} - \frac{1}{n}J_n. 
\end{equation*}
	
The transposition of matrix $A$ (or a vector) is denoted by $A^{\text T}$. For a square matrix $A$, $\Tr(A)$ refers to the summation of on-diagonal elements of $A$. The direct sum of any pair of matrices $A$ and $B$ is defined by
	\begin{equation*}
		A \oplus B ~:=~ \left[ \begin{array}{ccc}
			A & \mathbf 0 \\
			\mathbf 0& B  \end{array} \right].
	\end{equation*}
\begin{definition}[pseudo-inverse]
For $A \in \mathbb R^{n \times m}$, the Moore-Penrose pseudo-inverse of $A$ is defined by $A^{\dag} \in \mathbb R^{m \times n}$ satisfying all the following conditions
\begin{itemize}
	\item $AA^{\dag}A~=~A$,
	\item $A^{\dag}AA^{\dag}~=~A^{\dag}$,	
	\item $(AA^{\dag})^{\text T}~=~AA^{\dag}$,
	\item$(A^{\dag}A)^{\text T}~=~A^{\dag}A$.	
\end{itemize}
Note that for any matrix $A$, there is exactly one matrix  $A^{\dag}$ that satisfies all above conditions.
\end{definition}

{\it Graph Theory:} Throughout this paper, we assume that all graphs are finite, simple, and undirected. A weighted graph $\G$ is represented by a triple $\G = (\V_{\G}, \E_{\G},w_{\G})$, where $\V_{\G}$  is the set of nodes, $\E_{\G} \subseteq \big\{\{i,j\}\big| \,i,j \in \V_{\G}, ~i \neq j \big\}$ is the set of edges, and $w_{\G}: \E_{\G} \rightarrow \R_+$ is the weight function. An unweighted graph $\G$ is a graph with constant weight function $w_{\G}(e)=1$ for all $e \in \E_{\G}$. For each $i \in \V_{\G}$, the degree of node $i$ is defined by 
	\begin{equation*}
		d_i ~:= \sum_{e=\{i,j\} \in \E_{\G}} w_{\G}(e).
	\end{equation*}
%
The sum of all edge weights in graph $\G$ is denoted by $W(\G)$. The adjacency matrix
 $\adj = [a_{ij}]$ of graph $\mathcal{G}$ is defined by setting $a_{ij} = w_{\G}(e)$ if $e=\{i,j\} \in \E_{\G}$, otherwise $\mathsf{a}_{ij}=0$. The Laplacian matrix of $\mathcal G$ is defined by $L_{\G} := D_{\G} - \adj$, where $D_{\G}=\mathbf{diag}(d_{1},\ldots,d_{n})$ is a diagonal matrix. The eigenvalues of a Laplacian matrix $L_{\G}$ are indexed in ascending order $\lambda_1 \leq \lambda_2 \leq \cdots \leq \lambda_n $. If $\G$ is connoted, then $\lambda_1=0$. The class of all connected graphs with $n$ nodes is denoted by $\mathbb{G}_{n}$. 
The centering graph is  a complete graph with Laplacian matrix $M_n$ and is denoted by $\M_n$. 

	For comparison purposes throughout the paper, we consider the standard graphs in Table \ref{table} in several occasions. Every one of these graphs has its own comparable characteristics. For instance, among all graphs in $\mathbb{G}_{n}$ a complete graph has the maximum number of edges and a star graph has the maximum number of nodes of degree one. A path graph is a tree with minimum number of nodes of degree one. We refer to reference \cite{Bondy} for more details and discussions. A tree is a connected graph on $n$ nodes and with exactly $n-1$ edges. An unicyclic graph is a connected graph with exactly one cycle. A $d$-regular graph is a graph where all nodes have identical degree  $d$.

A subgraph $\mathcal{P}$ is a spanning subgraph of a graph $\G$ if it has the same node set as $\G$. An edge is called a cut-edge whose deletion increases the number of connected components. 

	\begin{table}[t]
		\centering
		\begin{tabular}{| c || c | c |}  
		\hline
		Standard Graph Families in $\mathbb{G}_{n}$ & Symbol\\
		\hline
		\hline 
		Complete   & $\mathcal K_n$ \\
		\hline
		Star  & $\mathcal S_n$\\
		\hline
		Cycle  & $\mathcal{C}_{n}$ \\
		\hline
		Path  & $\mathcal P_n$ \\
		\hline
		Bipartite & $\mathcal B_{n_1,n_2}$\\
		\hline
		Complete bipartite   & $\mathcal K _{n_1,n_2}$ \\
		\hline
		\end{tabular}
		\caption{\small{For comparison purposes throughout the paper, we consider the standard graphs in this table in several occasions. }}
		\label{table}
	\end{table}


\section{Linear Consensus Networks and their Performance Measures }
\label{sec:first-order}

	We consider two classes of linear consensus networks: first-order and second-order. The mathematical formulation is analogous in both cases, with the main difference being that the second-order models have two scalar states (position and velocity) locally for each subsystem contrary to a single scalar local state in the first-order model. These two classes of linear consensus networks have the following common canonical form
\begin{eqnarray}
		\NN(A;L_{\Q}):
		\begin{cases}
		~\dot \psi~=~ -A \psi~+~B \xi,\\
		~y~=~ C_{\mathcal{Q}} \psi,
		\end{cases}
		\label{consensus-general}
\end{eqnarray}
where  $\psi$ is the vector of state variables, $\xi$ is an exogenous uncorrelated white stochastic process with zero-mean and identity covariance matrix that can model random
forcing, $y$ is the performance output of the network, $A$ is the state matrix of the network, $C_{\mathcal{Q}}$ is the output matrix of the network and 
\begin{equation}
		L_{\mathcal{Q}}~=~C_{\mathcal{Q}}^{\text T}C_{\mathcal{Q}}.\label{C-L}
\end{equation}
The input matrix for the first-order consensus networks is $B=I$ and for the second-order consensus networks is $B=\big[\begin{array}{cc}\mathbf{0} & I\end{array}\big]^{\text T}$. 
	
\begin{definition}\label{output-graph}
A given graph $\mathcal{Q}=(\V_{\Q}, \E_{\mathcal{Q}},w_{\mathcal{Q}})$ 
is the output graph of a linear consensus network $\NN(A;L_{\Q})$ if $\Q$ admits $L_{\mathcal{Q}}$ in \eqref{C-L} as its Laplacian matrix. 
\end{definition}

In general, the output graph can be a disconnected graph. The output graphs help us to better understand how the specific choice of performance output will affect a given performance measure. In this paper, our primary goal is to discover existing relationships between a class of performance measures and the interconnection topology of the underlying graphs of the network. In the first step  toward this goal, we employ a class of performance measures that are defined using the performance outputs of linear consensus networks. 
 
\begin{definition}
Suppose that $\Q$ is an output graph of $\NN(A; L_{\Q})$. The steady state variance of the performance output of  the network $\NN(A; L_{\Q})$ is considered  as the performance measure
	\begin{equation}
	\RR(A; L_{\Q})~=~ \lim_{t \rightarrow \infty} \mathbb E \big[ y(t)^{\text T}y(t) \big]. \label{perf-meas}
	\end{equation}
\end{definition}

	In order to guarantee that the performance measure \eqref{perf-meas} is well-defined, marginally stable and unstable modes of $\NN(A; L_{\Q})$  should be unobservable from the performance output $y$. This can be true according to Definition \ref{output-graph} and the assumption that $\Q$ is an output graph for $\NN(A; L_{\Q})$.
\vspace{0.1cm}
\begin{assumption}\label{assump-output-graph}
For every linear consensus network $\NN(A; L_{\Q})$ considered in this paper, it is assumed that $L_{\Q}$ is a Laplacian matrix that corresponds to an output graph $\Q$. 
\end{assumption}	
\vspace{0.1cm}

The performance index \eqref{perf-meas} measures the performance of the network in the average. This is because \eqref{perf-meas} is indeed equivalent to the square of the $\mathcal H_2$--norm of the system from the exogenous disturbance input  to the performance output \cite{Bamieh12}. When there is no exogenous disturbance (noise) input, the steady state of $\NN(A;L_{\Q})$ in (\ref{consensus-general}) converges to the consensus state and the value of the performance measure becomes zero. In the following, we quantify performance measure \eqref{perf-meas} for each class of consensus networks.

\subsection{First-Order Consensus Networks}
The first class of consensus networks that we consider in this paper is the class of first-order consensus (FOC) networks whose dynamics are defined over  {\it coupling}  graphs $\mathcal G = (\V_{\G}, \E_{\G}, w_{\G})$ with $n$ nodes. For this class of networks, each node corresponds to a subsystem with a scalar state variable and the interconnection topology between these subsystems is defined by the coupling graph $\G$. The state of the entire network is represented by $x=[\begin{array}{cccc}x_{1} & x_{2} & \ldots & x_{n} \end{array}]^{\text T}$ where $x_{i}$ is the state variable of subsystem with index $i$ for  $i=1,\ldots,n$. For the class of FOC networks, the general consensus model \eqref{consensus-general} with the corresponding vector of state variable $\psi=x$ reduces to  
	\begin{eqnarray}
		\NN(L_{\G};L_{\Q}):
		\begin{cases}
		~\dot x~=~ -L_{\G}x~+~\xi,\\
		~y~=~ C_{\mathcal{Q}}x,
		\end{cases}
		\label{first-order-1}
	\end{eqnarray}
where  $A=L_{\G}$ is the Laplacian matrix of the coupling graph $\mathcal G$ and the input matrix is $B=I$.  

\vspace{0.1cm}
\begin{assumption}\label{assum-coupling-graph}
For every FOC network $\NN(L_{\G};L_{\Q})$ considered in this paper, the corresponding coupling graph $\G$ is assumed to be connected, i.e., $\G \in \mathbb{G}_n$.
\end{assumption}
\vspace{0.1cm}

	Based on Assumption \ref{assum-coupling-graph}, one can verify that  the state matrix of the network has exactly one marginally stable mode, which is unobservable from the performance output $y$. Because the output matrix of the network satisfies  the following property  
\begin{equation*}
C_{\mathcal{Q}}\mathbf 1=\mathbf 0.
\end{equation*} 
Therefore, the performance measure \eqref{perf-meas} is well-defined.

According to Assumption  \ref{assump-output-graph},   $L_{\mathcal{Q}}$  is the Laplacian matrix of the output graph $\Q$. Examples of admissible output matrices include incidence and centering matrices. In the special case where the output matrix is the centering matrix, i.e., $C_{\Q}~=~M_n$, it follows that
\begin{equation*}
		L_Q~=~C_{\Q}^{\text T}C_{\Q}~=~{M_n},
\end{equation*}
which implies that the output graph is a centering graph and we have $\Q=\M_n$.

\begin{theorem}
\label{1-thm}
For the FOC network (\ref{first-order-1}), the performance measure \eqref{perf-meas} can be quantified as   
	\begin{equation}
		\RR(L_{\G}; L_{\Q})~=~ \frac{1}{2}\mathbf{Tr}(L_{\mathcal Q}L_{\mathcal{G}}^{\dag}), \label{H-1}
	\end{equation}
where $L_{\mathcal{G}}^{\dag}$ is the Moore--Penrose pseudo inverse of $L_{\mathcal{G}}$. 
\end{theorem}
%
%
%

If the output graph is a centering graph, then the performance measure \eqref{H-1} reduces to 
	\begin{equation}
		\RR(L_{\G}; M_n)~=~\frac{1}{2}\Tr(L_{\G}^{\dag})~=~\frac{1}{2} \sum_{i=2}^{n}  \lambda_i^{-1},
		\label{eq-h}
	\end{equation}
where $\lambda_{i}$ for $i=2,\ldots,n$ are nonzero eigenvalues of $L_{\G}$ and $\lambda_1=0$ according to Assumption \ref{assum-coupling-graph}. The performance measure \eqref{H-1} relates to the concept of coherence in consensus networks  and the expected dispersion of the state of the system in steady state \cite{Young10, Bamieh12}. The performance measure \eqref{H-1} has also close connections to the notion of the total effective resistance  of the coupling graph $\G$ as follows
	\begin{equation}
		\RR(L_{\G}; M_n)~=~\frac{\mathbf{r}_{\textrm{total}}}{2n},
		\label{eq-r}
	\end{equation}
where the total effective resistance  of  $\G$ is given by 
\cite{Bamieh12, Barooah}
	\begin{equation*}
		\mathbf{r}_{\textrm{total}}~=~n\sum_{i=2}^{n} \lambda_i^{-1}.
	\end{equation*}

\subsection{Second-Order Consensus Networks}\label{sec:second-order}
The second class of consensus networks that we consider in this paper is the class of second-order linear consensus (SOC) networks. For a given SOC network, each node corresponds to a subsystem with two scalar state variables $x_i$ and $v_i$ for $i=1,\ldots,n$. The state of the entire network is obtained by concatenating the state vectors $x=[\begin{array}{cccc}x_{1} & x_{2} & \ldots & x_{n} \end{array}]^{\text T}$ and $v=[\begin{array}{cccc}v_{1} & v_{2} & \ldots & v_{n} \end{array}]^{\text T}$. For the class of SOC networks, the general consensus model \eqref{consensus-general} with the corresponding vector of state variable $\psi=[\begin{array}{cc}x^{\text T} & v^{\text T}\end{array}]^{\text T}$ takes the following form
\begin{equation}
\NN(A;L_{\Q}):
\begin{cases}
		\left[ \begin{array}{c}
		\dot x \\
		\dot v \end{array} \right] =\left[ \begin{array}{ccc}
		\mathbf 0& I  \\
		F & G  \end{array} \right] \left[ \begin{array}{c}
		x\\
		v \end{array} \right]+\left[ \begin{array}{c}
		\mathbf 0\\
		I \end{array} \right] \xi  \\
		\\
\hspace{0.75cm}		y = C_{\Q}
		\left[ \begin{array}{c}
 		x \\
 		v \end{array} \right] 
		\end{cases}		
		\label{second-order}
\end{equation}
where the state and input matrices are replaced  by 
	\begin{equation} 
		A~=~-\left[ \begin{array}{ccc}
		\mathbf 0& I  \\
		F & G  \end{array} \right], ~~~B=\left[ \begin{array}{c}
		\mathbf 0\\
		I \end{array} \right]. \label{A-matrix}
	\end{equation}
It is assumed that $F$ and $G$ are some stabilizing static linear feedback gains and their  specific structure depend on the types of sensor measurements used to close the feedback loop \cite{Bamieh12}. We impose the following structural constraint on the output matrix of the network  
	\begin{equation}
		C_{\Q}~=~C_{\Q_x} \oplus C_{\Q_v},
		\label{output-struct}
	\end{equation}
where $C_{\Q_x} \in \R^{n \times n}$ is the position output matrix and $C_{\Q_v} \in \R^{n \times n}$ is the velocity output matrix. According to Assumption  \ref{assump-output-graph},  $L_{\Q}$ is a Laplacian matrix with the following decomposition 
	\begin{equation*}
		L_{\Q}~=~C_{\Q_x}^{\text T}C_{\Q_x} \oplus C_{\Q_{v}}^{\text T}C_{\Q_{v}}~=~L_{\Q_x} \oplus L_{\Q_v}.
	\end{equation*}
This implies that $L_{\Q_{x}}$ and $L_{\Q_{v}}$ are also Laplacian matrices. Thus, we associate two output graphs to a given SOC network and name $\Q_{x}$
as the position output graph with Laplacian  $L_{\Q_{x}}$, and $\Q_{v}$ as the  velocity output graph with Laplacian $L_{\Q_{v}}$. Based on our notations, the performance measure of \eqref{second-order} can be decomposed with respect to output graphs as follows 
	\begin{equation*}
		\RR (A;L_{\Q})=\RR (A; L_{\mathcal{Q}_{x}} \oplus \mathbf 0)+\RR (A; \mathbf 0 \oplus L_{\mathcal{Q}_{v}}).
	\end{equation*}
In this paper, our primary focus is on two types of SOC networks that each is defined over a given coupling graph $\mathcal G = (\V_{\G}, \E_{\G}, w_{\G})$ with $n$ nodes and Laplacian matrix $L_{\G}$. For each type, the state matrix $A$ in \eqref{A-matrix} is defined as follows
\begin{eqnarray}
 & \textrm{\bf Type 1:} & \hspace{2cm} A^{(1)}_{\G} := \left[ \begin{array}{ccc}
		\mathbf 0& -I  \\
		L_{\G} & \beta I  \end{array} \right],
		\label{type-1} \\
& \textrm{\bf Type 2:} & \hspace{2cm}  A^{(2)}_{\G} := 
		\left[ \begin{array}{ccc}
		\mathbf 0& -I  \\
		L_{\G} & \beta L_{\G}  \end{array} \right], \label{type-2}
\end{eqnarray}
for some design parameter $\beta > 0$. These consensus models are previously studied in \cite{Bamieh12}. In Section \ref{sec:applications}, we will discuss two real-world applications and show that a linearized model of a power network can be cast as a Type 1 SOC network and an abstract model of controlled vehicles in a formation can be written in the form of a Type 2 SOC network. 

\vspace{0.1cm}
\begin{assumption}\label{assum-coupling-graph-SOC}
For all Type 1 and Type 2  SOC networks   considered in this paper, it is assumed that 
the corresponding coupling graphs are  connected, i.e., $\G \in \mathbb{G}_n$. 
\end{assumption}
\vspace{0.1cm}

	From this assumption, it follows that the state matrices $A^{(1)}_{\G}$ has a simple zero eigenvalue and all the other eigenvalues have positive real parts, and $A^{(2)}_{\G}$ has a zero eigenvalue of multiplicity two and all the other eigenvalues have positive real parts (cf. \cite{Yu2010}).
According to assumption \eqref{output-struct}, the performance measure \eqref{perf-meas} is well-defined for these two types of SOC networks. 

	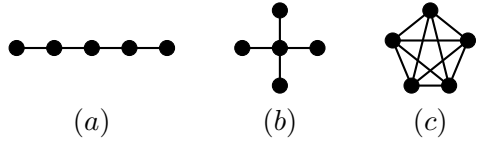
\begin{figure}[t]
		\centering
            \begin{tikzpicture}
            \draw [fill] (-4,0) circle [radius=0.1];
            \draw [fill] (-3.5,0) circle [radius=0.1];
            \draw [fill] (-3,0) circle [radius=0.1];
            \draw [fill] (-2.5,0) circle [radius=0.1];
            \draw [fill] (-2,0) circle [radius=0.1];
            \draw [ thick] (-4,0) -- (-2,0);
            \draw [fill] (-.5,0) circle [radius=0.1];
            \draw [fill] (-1,0) circle [radius=0.1];
            \draw [fill] (0,0) circle [radius=0.1];
            \draw [fill] (-.5,0.5) circle [radius=0.1];
            \draw [fill] (-.5,-0.5) circle [radius=0.1];
            \draw [ thick] (-1,0) -- (0,0);
            \draw [ thick] (-.5,-.5) -- (-.5,.5);
            \draw [fill] (1,0.1) circle [radius=0.1];
            \draw [fill] (2,0.1) circle [radius=0.1];
            \draw [fill] (1.5,0.5) circle [radius=0.1];
            \draw [fill] (1.75,-0.5) circle [radius=0.1];
            \draw [fill] (1.25,-0.5) circle [radius=0.1];
            \draw [ thick] (1,0.1) -- (2,0.1);
            \draw [ thick] (1,0.1) -- (1.5,0.5);
            \draw [ thick] (1,0.1) -- (1.75,-0.5);
            \draw [ thick] (1,0.1) -- (1.25,-0.5);
            \draw [ thick] (2,0.1) -- (1.5,0.5);
            \draw [ thick] (2,0.1)-- (1.75,-0.5);
            \draw [ thick] (2,0.1)-- (1.25,-0.5);
            \draw [ thick] (1.5,0.5) -- (1.75,-0.5);
            \draw [ thick] (1.5,0.5) -- (1.25,-0.5);
            \draw [ thick] (1.75,-0.5)-- (1.25,-0.5);
            \node[] at (-3,-1) {$(a)$};
            \node[] at (-.5,-1) {$(b)$};
            \node[] at (1.5,-1) {$(c)$};
            \end{tikzpicture}
            \caption{ This figure illustrates the results of  Theorems \ref{max-min-thm} and \ref{tree-thm} for  the following extreme cases.  The performance measure \eqref{rho-ss} is (a) maximal for $\mathcal P_{5}$ among all graphs as well as among all trees in $\mathbb{G}_5$, (b) minimal for $\mathcal S_5$ among all trees  in $\mathbb{G}_5$, and (c) minimal for $\mathcal K_5$ among all graphs in $\mathbb{G}_5$.} \label{fig_212}
\end{figure}

\begin{theorem}
\label{2-thm}
For type 1 and 2 SOC networks with state matrices \eqref{type-1} and \eqref{type-2}, the performance measure \eqref{perf-meas} can be quantified as   
	\begin{eqnarray}
		\RR (A^{(1)}_{\G}; L_{\Q_x} \oplus \mathbf 0) & = & \frac{1}{2\beta}\mathbf{Tr}(L_{\Q_{x}}L_{\G}^{\dag}),\label{position-1} \\
		\RR (A^{(1)}_{\G};\mathbf 0 \oplus  L_{\Q_v}) & = & \frac{1}{2\beta}\mathbf{Tr}(L_{\Q_{v}}), \label{velocity-1} \\
		\RR (A^{(2)}_{\G}; L_{\Q_x} \oplus \mathbf 0) & = & \frac{1}{2\beta}\mathbf{Tr}(L_{\Q_{x}}(L_{\G}^{\dag})^2),
		\label{position-2} \\
		\RR (A^{(2)}_{\G};\mathbf 0 \oplus  L_{\Q_v}) & =&\frac{1}{2\beta}\mathbf{Tr}(L_{\Q_{v}}L_{\G}^{\dag}),
		\label{velocity-2}
	\end{eqnarray}
where $L_{\mathcal{G}}^{\dag}$ is the Moore--Penrose pseudo inverse of $L_{\mathcal{G}}$. 
\end{theorem}

The steps involved for calculating performance measures  \eqref{position-1} and \eqref{velocity-2} are identical to those of \eqref{H-1} for the FOC networks. The value of performance measure \eqref{velocity-1} only depends on the velocity output graph $\Q_v$ and is equal to $\frac{W(\Q_v)}{\beta}$, where $W(\Q_v)$ is the sum of all edge weights of $\Q_v$. For performance measure \eqref{position-2}, if the corresponding position output graph $\Q_x$ is a centering graph, then the performance measure \eqref{position-2} reduces to 
	\begin{equation}
		\RR (A^{(2)}_{\G}; M_n \oplus \mathbf 0) ~ = ~ \frac{1}{2\beta}\mathbf{Tr}((L_{\G}^{\dag})^2)~=~\frac{1}{2\beta} \sum_{i=2}^{n} \lambda_i^{-2},
		\label{perf-meas-soc}
	\end{equation}
where $\lambda_{i}$ for $i=2,\ldots,n$ are nonzero eigenvalues of $L_{\G}$ and $\lambda_1=0$ according to Assumption \ref{assum-coupling-graph-SOC}. We end this section by summarizing that performance analysis of FOC and SOC networks involves calculating the following spectral-zeta functions 
	\begin{equation*}
		\zeta_{\G}(p) \hspace{0.05cm} = \hspace{0.05cm} \sum_{i=2}^{n}\hspace{0.05cm} \lambda_i^{-p},
	\end{equation*}
for $p \in \{1,2\}$. 

\section{Fundamental Limits on the Performance Measure }\label{sec-IV}

We evaluate the performance of FOC networks \eqref{first-order-1} with respect to  the centering output graph by considering the following performance measure 
	\begin{equation}
\RR (L_{\G}; M_n) ~=~ \frac{1}{2} \sum_{i=2}^n \lambda_i^{-1}.  \label{rho-ss}
	\end{equation}
This consideration also covers analysis of performance measure \eqref{position-1} with the   centering position output graph as well as performance measure \eqref{velocity-2} with the   centering velocity output graph. Therefore, among all four performance measures \eqref{position-1}-\eqref{velocity-2} for SOC networks, it only remains  to treat performance measure \eqref{perf-meas-soc} separately. In this section, several scenarios are investigated in order to reveal the fundamental role of the coupling graphs of FOC and SOC networks on emergence of fundamental limits on these performance measures. 

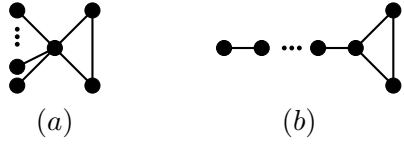
\begin{figure}[t]
		\centering
		\begin{tikzpicture}
        	\draw [ thick] (-1,-.25) -- (-.5,0);
        	\draw [ thick] (-1,.5) -- (-.5,0);
        	\draw [ thick] (-1,-.5) -- (-.5,0);
        	\draw [ thick] (0,-.5) -- (0,.5);
        	\draw [ thick] (-.5,0) -- (0,.5);
        	\draw [ thick] (0,-.5) -- (-.5,0);
        	\draw [fill] (-.5,0) circle [radius=0.1];
        	\draw [fill] (0,.5) circle [radius=0.1];
        	\draw [fill] (0,-.5) circle [radius=0.1];
        	\draw [fill] (-1,.5) circle [radius=0.1];
        	\draw [fill] (-1,-.25) circle [radius=0.1];
        	\draw [fill] (-1,-.5) circle [radius=0.1];
		\draw [fill] (-1, .25) circle [radius=0.1/4];
		\draw [fill] (-1, .15) circle [radius=0.1/4];
		\draw [fill] (-1, 0.05) circle [radius=0.1/4];
	
        	\draw [ thick] (-1+4,0) -- (-.5+4,0);
		\draw [ thick] (-2.25+4,0) -- (-1.75+4,0);
        	\draw [ thick] (0+4,-.5) -- (0+4,.5);
        	\draw [ thick] (-.5+4,0) -- (0+4,.5);
        	\draw [ thick] (0+4,-.5) -- (-.5+4,0);
        	\draw [fill] (-.5+4,0) circle [radius=0.1];
        	\draw [fill] (0+4,.5) circle [radius=0.1];
        	\draw [fill] (0+4,-.5) circle [radius=0.1];
        	\draw [fill] (-1+4,0) circle [radius=0.1];
		\draw [fill] (-1.75+4,0) circle [radius=0.1];
		\draw [fill] (-2.25+4,0) circle [radius=0.1];
		\draw [fill] (-1.75+4+.3,0) circle [radius=0.1/4];
		\draw [fill] (-1.75+4+.4,0) circle [radius=0.1/4];
		\draw [fill] (-1.75+4+.5,0) circle [radius=0.1/4];
		\node[] at (-.5,-1) {$(a)$};
           \node[] at (-2.5/2+4,-1) {$(b)$};
        	\end{tikzpicture}
        	\caption{{\small The unicyclic graphs that achieve the lower and upper bounds in Theorem \ref{unicyclic-thm}: (a) $\mathcal{G}=\mathcal{S}(\K_3; \mathcal K_1, \cdots ,\mathcal K_1)$, and (b) $\mathcal{P}(\K_3; \mathcal K_1, \cdots,\mathcal K_1)$. }}
         \label{fig_67}
	\end{figure}
	
\subsection{Universal Bounds and Scaling Laws}
\label{universal-sub}

The following result presents universal lower and upper bounds for the best and worst achievable values for performance measure \eqref{rho-ss} among all FOC networks with arbitrary unweighted coupling graphs.

\begin{theorem} \label{max-min-thm}
For a given FOC network with an unweighted coupling graph $\G \in \mathbb G_n$,  the  performance measure \eqref{rho-ss} is universally bounded by%
	\begin{equation}
		\frac{1}{2} -\frac{1}{2n} \hs \leq \hs  \RR (L_{\G}; M_n) \hs \leq \hs \frac{n^2-1}{12}.
		\label{eq-the-3}
	\end{equation}
Furthermore, the lower bound is achieved if and only if $\mathcal G=\mathcal{K}_{n}$, and the upper bound is reached if and only if $\mathcal G=\mathcal{P}_{n}$.
\end{theorem}

The bounds in inequalities \eqref{eq-the-3}  are conservative as they only depend on the network size and nothing specific is known about the topology of the coupling graph of the network.  These bounds can be tightened if we consider more specific subclasses of graphs. In the following three theorems, we improve the bounds in Theorem  \ref{max-min-thm} for three important classes of graphs: tree, unicyclic, and bipartite.

\begin{theorem}
\label{tree-thm}
For a given FOC network with an unweighted tree coupling graph $\mathcal{T} \in \mathbb{G}_n$ with $n \geq 5$, the  performance measure \eqref{rho-ss} is bounded by
	\begin{equation}
		\frac{(n-1)^2}{2n} \hs \leq \hs \mathbf  \RR (L_{\mathcal T}; M_n) \hs \leq \hs  \frac{n^2-1}{12}.
		\label{tree-inq}
	\end{equation}
Moreover, the lower bound is achieved if and only if  $\mathcal{T}=\mathcal{S}_n$, and the upper bound is achieved if and only if  $\mathcal{T}=\mathcal{P}_n$.
\end{theorem}

The lower bound in \eqref{tree-inq} is tight as if the value of the performance measure is strictly less than $\frac{(n-1)^2}{2n}$, then the unweighted coupling graph of the network must have at least one cycle. The next result quantifies  the effects of adding exactly one cycle to coupling graph of a FOC network whose coupling graph is a tree. 

\begin{theorem}\label{unicyclic-thm}
For a given FOC network with an unweighted unicyclic coupling graph in $\mathbb G_n$ with $n \geq 13$, the  performance measure \eqref{rho-ss} is bounded by
	\begin{equation}
		\frac{(n-1)^2}{2n} - \frac{1}{3} \hs \leq \hs  \RR (L_{\G};M_n) \hs \leq \hs \frac{n^2-1}{12}+\frac{3}{2n}-1.
		\label{uni-inq}
	\end{equation}
Moreover, the lower bound is achieved if and only if  $\mathcal{G}=\mathcal{S}(\K_3; \mathcal K_1, \cdots ,\mathcal K_1)$, which is a star-like graph that is formed by replacing the center of $\mathcal{S}_n$ by a clique $\mathcal K_3$, and the upper bound is achieved if and only if  $\mathcal{G}=\mathcal{P}(\K_3; \mathcal K_1, \cdots,\mathcal K_1)$, which is a path-like graph that is formed by replacing one of the end nodes of $\mathcal{P}_n$ by a clique $\mathcal K_3$.
\end{theorem}

The lower and upper bounds in \eqref{uni-inq} are tight, in the sense that if the value of the performance measure for a FOC network does not satisfy \eqref{uni-inq}, then the coupling graph of this network is either a tree (with no cycle) or has at least two cycles.

The following theorem investigates the performance of a class of FOC networks defined over bipartite graphs. In this case the network consists of two disjoint sets of nodes and the states of one set depend on the states of the other set and vice versa. Bipartite graphs appear in many applications, for instance when we model a network of electricity sellers and buyers \cite[Ch.12]{networks-book}, a power network \cite[Sec. 2]{Kraning2014}, and a network of leaders and followers agents where leaders only influenced by their followers and vice versa.

	\begin{table}[t]
       	\centering
       	\begin{tabular}{ |c || c | c|  }  
       	\hline
       	Unweighted Coupling Graphs  &Lower Bound & Upper Bound\\
       	\hline
       	\hline 
       	Arbitrary & $\frac{1}{2}-\frac{1}{2n}$ & $\frac{n^2-1}{12}$ \\
       	\hline
       	Tree &$\frac{(n-1)^2}{2n}$ & $\frac{n^2-1}{12}$ \\
       	\hline
       	Unicyclic &$\frac{(n-1)^2}{2n} - \frac{1}{3}  $ & $\frac{n^2-1}{12}+\frac{3}{2n}-1$ \\
       	\hline
		Bipartite &$ 1-\frac{\lfloor \frac{n}{2} \rfloor }{n\lceil \frac{n}{2} \rceil }$ & $\frac{n^2-1}{12}$ \\
       	\hline
        	\end{tabular}
        	\caption{\small The universal bounds on the performance measure \eqref{rho-ss} for  FOC networks with unweighted coupling graphs in $\mathbb G_n$.}
        	\label{table1}
	\end{table}

\begin{theorem}\label{bipartite-thm} 
For a given FOC network with an unweighted  bipartite coupling graph $\mathcal{B}_{n_1,n_2} \in \mathbb G_n$ with $n_1 + n_2 =n$, the  performance measure \eqref{rho-ss} is bounded by
	\begin{equation}
		 1-\frac{\lfloor \frac{n}{2} \rfloor }{n\lceil \frac{n}{2} \rceil } \hs \leq \hs  \RR(L_{\mathcal{B}};M_n) \hs \leq \hs \frac{n^2-1}{12}.
		\nonumber
	\end{equation}
Furthermore, the lower bound is achieved if and only if   $\mathcal{B}_{n_1,n_2} = \mathcal K_{\lfloor \frac{n}{2} \rfloor, \lceil \frac{n}{2} \rceil}$, and the upper bound is achieved if and only if  $\mathcal{B}_{n_1,n_2} =\mathcal P_n$, where $\lfloor \hspace{0.03cm}.\hspace{0.03cm} \rfloor$ and $\lceil \hspace{0.03cm}.\hspace{0.03cm} \rceil$ are the floor and ceiling operators, respectively. 
\end{theorem}

The lower bound in Theorem \ref{bipartite-thm} is tight. This is because if the value of the performance measure is strictly less than $1-\frac{\lfloor \frac{n}{2} \rfloor }{n\lceil \frac{n}{2} \rceil }$ for a given FOC network with an unweighted coupling graph, then the coupling graph of the network cannot be a complete bipartite graph.

In the following result, we quantify universal lower and upper bounds for the best and worst achievable values for performance measure \eqref{perf-meas-soc} among all Type 2 SOC networks with arbitrary unweighted coupling graphs.  

\begin{theorem}
\label{general-thm-type-2}
For a given Type 2 SOC network with an unweighted coupling graph $\mathcal{G} \in \mathbb G_n$, the  performance measure (\ref{perf-meas-soc}) is bounded by
	\begin{equation}
		\frac{1}{\beta}\left( \frac{1}{2n}-\frac{1}{2n^2}\right) \hs \leq  \hs \RR (A^{(2)}_{\G}; M_n \oplus \mathbf 0)  \hs < \hs  \Xi(n),
		\label{SOC-general-inq}
	\end{equation}
where
	\begin{equation*}
		\Xi(n) = \frac{(n^2-1)^2}{72\beta}-\frac{(n-1)(n-2)^2}{2n\beta}.
	\end{equation*}
Moreover, the lower bound is achieved if and only if  $\mathcal{G}=\mathcal K_n$.
\end{theorem}

According to our current analysis, the upper bound in Theorem \ref{general-thm-type-2} is not tight as the performance measure is strictly less than the hard limit function $\Xi(n)$. By comparing this result to the result of Theorem \ref{max-min-thm}, we observe that the performance measure \eqref{perf-meas-soc} does not attain its minimum value for SOC networks with path coupling graphs. One can verify that the performance measure of a Type 2 SOC network with path coupling graph is given by 
\begin{equation}
\RR (A^{(2)}_{\mathcal P_n}; M_n \oplus \mathbf 0) \hs = \hs \frac{(n^2-1)^2}{72\beta}-\frac{1}{n\beta}\binom{n+2}{5}. 
\label{SOC-path}
\end{equation}
In the following, we obtain less conservative lower bound by limiting our attention to the class of tree graphs.

\begin{theorem}
\label{tree-thm-type-2}
For a given Type 2 SOC network with an unweighted tree coupling graph $\mathcal{T} \in \mathbb G_n$ with $n \geq 5$, the  performance measure (\ref{perf-meas-soc}) is bounded by
	\begin{equation}
		\frac{1}{2\beta}+\frac{(n-2)^3}{2\beta(2n-3)^2} \hs <  \hs \RR (A^{(2)}_{\mathcal T}; M_n \oplus \mathbf 0)  \hs < \hs  \Xi(n).
		\label{SOC-tree-inq}
	\end{equation}
\end{theorem}

Based on our current proof methods, both lower and upper bounds in Theorem \ref{tree-thm-type-2} are not tight. In contrary to Theorem \ref{tree-thm}, the lower bound for the performance measure \eqref{perf-meas-soc} among all trees does not happen for a Type 2 SOC network  with  star coupling graph, which is given by 
\begin{equation} 
\RR (A^{(2)}_{\mathcal S_n}; M_n \oplus \mathbf 0) \hs = \hs \frac{1}{\beta}\left (\frac{n}{2}+\frac{1}{2n^2}-1\right). \label{SOC-star}
\end{equation}
Nevertheless, we conjecture that \eqref{SOC-path} and \eqref{SOC-star} provide tight upper and lower bounds for \eqref{SOC-tree-inq}, and \eqref{SOC-path} provides a tight upper bound for \eqref{SOC-general-inq} as well.

	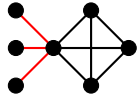
\begin{figure}[t]
		\centering
		\begin{tikzpicture}
        \draw [ thick] (-.5,0) -- (.5,0);
        \draw [ thick, red] (-1,0) -- (-.5,0);
        \draw [ thick, red] (-1,.5) -- (-.5,0);
        \draw [ thick, red ] (-1,-.5) -- (-.5,0);
        \draw [ thick] (0,-.5) -- (0,.5);
        \draw [ thick] (-.5,0) -- (0,.5);
        \draw [ thick] (0,.5) -- (.5,0);
        \draw [ thick] (0,-.5) -- (.5,0);
        \draw [ thick] (0,-.5) -- (-.5,0);
        \draw [fill] (-.5,0) circle [radius=0.1];
        \draw [fill] (.5,0) circle [radius=0.1];
        \draw [fill] (0,.5) circle [radius=0.1];
        \draw [fill] (0,-.5) circle [radius=0.1];
        \draw [fill] (-1,.5) circle [radius=0.1];
        \draw [fill] (-1,0) circle [radius=0.1];
        \draw [fill] (-1,-.5) circle [radius=0.1];
        \end{tikzpicture}
        \caption{{\small A schematic graph of $\mathcal S(\mathcal K_4; \mathcal K_1,\mathcal K_1,\mathcal K_1)$ that has the minimal value of performance measure among all graphs in $\mathbb{G}_7$ with exactly $3$ cut edges (highlighted by red color).}}
         \label{fig_22}
	\end{figure}

 We summarize our results in this part by observing that  according to Theorems \ref{max-min-thm} and \ref{bipartite-thm} the performance measure of large-scale FOC networks with fairly  dense coupling graphs scale at least constantly with the network size as follows
	\begin{equation*}
		\RR (L_{\G}; M_n) \hs = \hs \Omega(1), 
	\end{equation*} 
and according to Theorems \ref{tree-thm} and \ref{unicyclic-thm} for large-scale FOC networks with fairly sparse coupling graphs, the performance measure scales at least linearly with the network size as follows
	\begin{equation*}
		\RR (L_{\G}; M_n) \hs = \hs \Omega(n).  
	\end{equation*} 
In the extreme cases, the performance measure of FOC networks scale at most quadratically with the network size as follows 
	\begin{equation*}
		\RR (L_{\G}; M_n) \hs = \hs \mathcal O(n^2). \label{scaling-2}
	\end{equation*} 
Based on Theorem \ref{general-thm-type-2}, the performance measure of large-scale Type 2 SOC networks with fairly dense coupling graphs decay  of the order of  	
	\begin{equation*}
		\RR (A^{(2)}_{\G}; M_n \oplus \mathbf 0) = \Omega (n^{-1}),	
	\end{equation*}
and based on Theorem \ref{tree-thm-type-2}, the performance measure of large-scale Type 2 SOC networks with fairly sparse coupling graphs scale at least linearly with the network size as follows
	\begin{equation*}
		\RR (A^{(2)}_{\G}; M_n \oplus \mathbf 0) = \Omega(n). 
	\end{equation*}
For the class of Type 2 SOC networks, the performance measure scales polynomially with the network size as follows 
	\begin{equation*}
		\RR (A^{(2)}_{\G}; M_n \oplus \mathbf 0) = \mathcal O(n^4). 
	\end{equation*}

%

	\begin{table}[t]
        \centering
        \begin{tabular}{ |c || c |  }  
        \hline
        Known Graph Specification & Lower Bound  \\
        \hline
        \hline 
        $\kappa(\G)$ & $\frac{1}{2n}+\frac{\kappa(\G)+1}{2}-\frac{1}{n-\kappa(\G)}$   \\
        \hline
        $\{d_i\}_{i=1}^n$&$-\frac{1}{2n}+\frac{n-1}{2n} \sum_{i=1}^{n} \frac{1}{d_i}$   \\
        \hline
        $\sigma(\G)$&$\frac{n-1}{2 \sigma(\G)}$  \\
        \hline
        $\mathfrak{T}(\mathcal G)$ &$\frac{n-1}{2\hspace{0.04cm} \sqrt[n-1]{n\mathfrak{T}(\mathcal G)}}$  \\
        \hline
        \end{tabular}
        \caption{\small A summary of our results that lists our combinatorial lower bounds on the best achievable values of performance measure for FOC networks with unweighted coupling graphs in  $\mathbb G_n$.}
        \label{table2}
	\end{table}

\subsection{Bound Calculations via Exploiting Structure of Coupling Graphs}
\label{tailored-sub}

In the previous subsection, we derived universal lower and upper bounds for performance measures of FOC and SOC networks with unweighted coupling graphs. Our results are summarized in Table \ref{table1}. These bounds are only functions of the network size. In this subsection, we incorporate additional known graph specifications in calculating lower and upper bounds for performance measures. We consider five important graph specifications and  expand our analysis through FOC and SOC networks with weighted and unweighted coupling graphs. Summaries of our results are listed in Tables \ref{table2} and \ref{table3}. 

\subsubsection{Graph diameter and number of edges}  The diameter of a graph is defined as the largest distance between every pair of nodes in that graph. 

\begin{theorem} \label{diam-thm}
For a given FOC network with an arbitrary  unweighted coupling graph $\G \in \mathbb G_n$, the performance measure  \eqref{rho-ss} is bounded by
	\begin{equation}
 		\mathfrak{L}_{\G} \hs \leq  \hs \RR (L_{\G}; M_n) \hs \leq \hs  \mathfrak{U}_{\G}, \label{gen-low-up-bound}
	\end{equation}
where
	\begin{eqnarray*}
		\mathfrak{L}_{\G}&=&\frac{(n-1)^{2}}{4 m} \label{gen-low}\\
		\mathfrak{U}_{\G}&=& \frac{1}{2n} \left [ n-1 +\left [\binom{n}{2}-m \right]\mathbf{diam}(\mathcal G) \right ]
	\end{eqnarray*}
where $\mathbf{diam}(\mathcal G)$ is the diameter and $m$ is the number of edges of $\G$. If $\G= \mathcal K_n$, then the lower and upper bounds in \eqref{gen-low-up-bound} coincide and  
	\begin{equation*}
		\RR (L_{\mathcal K_n};M_n)~=~\frac{n-1}{2n}.
	\end{equation*}  
Moreover, a star graph $\mathcal S_n$ achieves the upper bound in \eqref{gen-low-up-bound}. 
\end{theorem}
%

\subsubsection{Total weight sum} The sum of all edge weights in a weighted graph $\G$ is defined by  
\[ W(\G) :=\sum_{e \in \E_{\G}}w_{\G}(e).\]

\begin{theorem}\label{trace-thm}
For a given FOC network with an arbitrary weighted coupling graph $\G \in \mathbb G_n$, the performance measure  \eqref{rho-ss} is bounded from below by
	\begin{equation}
		\RR (L_{\G};M_n) \hs \geq \hs \frac{(n-1)^2}{4 W(\G)}.
		\label{propo}
	\end{equation}
\end{theorem}

	\begin{table}[t]
        \centering
        \begin{tabular}{ |c || c | c|  }  
        \hline
         Known Graph Specification &   Lower Bound    &  Upper Bound  \\
        \hline
        \hline 
        $W(\G)$ & $\frac{(n-1)^2}{4W(\G)}$ & $\infty$\\
   \hline      
                $\mathfrak{T}(\mathcal G)$ &$\frac{n-1}{2\hspace{0.04cm} \sqrt[n-1]{n\mathfrak{T}(\mathcal G)}}$  & $\infty$ \\      
        
        \hline
        $\{d_i\}_{i=1}^n$ &$ \displaystyle \Delta(\G) $ & $\infty$ \\
        \hline
        \end{tabular}
        \caption{\small A summary of our results that lists our combinatorial lower bounds on the best achievable values of performance measure for FOC networks with weighted coupling graphs in  $\mathbb G_n$.}
        \label{table3}
	\end{table}

\subsubsection{Number of spanning trees}
A spanning subgraph of $\G$ is called a spanning tree if it is also a tree. The weighted number of spanning trees of a connected graph $\G = (\V_{\G}, \E_{\G},w_{\G})$ 
is defined by
	\begin{equation} 
		\mathfrak{T}({\mathcal G}) \hs := \hs  \sum_{\mathcal T}\prod_{e \in \E_{\mathcal T}} w_{\G}(e),  \label{no-span-trees}
	\end{equation}
where the summation runs over all spanning trees $\mathcal T$ of $\mathcal G$.  For unweighted graphs, the total number of spanning trees of a connected graph is an invariant graph specification. 

\begin{theorem}
\label{spanning-thm}
For a given FOC network with an arbitrary weighted coupling graph $\G \in \mathbb G_n$, the performance measure  \eqref{rho-ss} is bounded from below by
	\begin{equation}
		\RR (L_{\G};M_n) \hs \geq \hs \frac{n-1}{2\hspace{0.04cm} \sqrt[n-1]{n\mathfrak{T}(\mathcal G)}},				\label{ineq-spanning-tree}
	\end{equation}
where $\mathfrak{T}(\mathcal G)$ is the  number of spanning trees of $\G$ defined by \eqref{no-span-trees}. 
\end{theorem}

The result of this theorem holds for general weighted connected graphs. However, for some particular classes of unweighted connected graphs, the total number of spanning trees can be calculated explicitly as a function of $n$. For example, for an unweighted complete graph $\mathcal K_n$ the total number of spanning trees is  $\mathfrak{T}(\mathcal G)=n^{n-2}$. In fact,   the lower bound in \eqref{ineq-spanning-tree} is tight for weighted graphs and it can be achieved by complete graphs. Nonetheless, our analysis shows that the proposed lower bound in \eqref{ineq-spanning-tree} is not tight for the class of unweighted tree, cycle, and complete bipartite graphs. As we discussed earlier, our results in Subsection \ref{universal-sub} are tight for these classes of graphs. 
  





\subsubsection{Number of cut edges}  An edge $e$ is called a cut edge of $\G$ if removing $e$ from $\G$ results in more connected components than $\G$.  The total number of cut edges in $\G$ is denoted by $\kappa(\G)$.

\begin{theorem}\label{cut-thm}
For a given FOC network with an arbitrary unweighted coupling graph $\G \in \mathbb G_n$ that has $\kappa(\G)$ cut edges, the performance measure  \eqref{rho-ss} is bounded from below by
	\begin{equation}\label{lower-cut-edge}
		\RR (L_{\mathcal G}; M_n) \hs \geq \hs \frac{1}{2n}+\frac{\kappa(\G)+1}{2}-\frac{1}{n-\kappa(\G)}.
	\end{equation}
%
The equality holds if and only if $\mathcal G = \mathcal S(\mathcal K_{n-\kappa(\G)}; \mathcal K_1, \cdots , \mathcal K_1)$, \ie $\mathcal G$ is a star graph that is formed by replacing the center of the star with a clique $\mathcal K_{n-\kappa(\G)}$.
\end{theorem}

For a given graph in $\mathcal{G}_n$, the number of cut edges satisfies  
\[ 0 \hs \leq \hs \kappa(\G) \hs \leq \hs n-1,\] 
where a tree with $n-1$ cut edges has the maximum  and a complete graph with zero cut edge has the minimum number of cut edges among all graphs in $\mathcal{G}_n$. A simple calculation reveals that the lower bound in \eqref{lower-cut-edge} gains its maximum value for tree and its minimum value for complete coupling graphs. This asserts that the lower bound in \eqref{lower-cut-edge}  is tight according to the results of Theorems  \ref{max-min-thm} and \ref{tree-thm}.

\subsubsection{Degree sequence}  

\begin{theorem}\label{degree-2-thm}
For a given FOC network with an arbitrary  weighted  coupling graph $\G \in \mathbb G_n$ and degree sequence $\big\{ d_i\big\}_{i=1}^n$, the  performance measure \eqref{rho-ss} is bounded from below by
\begin{equation}
\RR (L_{\G};M_n) \hs \geq \hs \Delta(\G), \label{lower-delta}
\end{equation}
where 
\begin{equation}
\Delta(\G):= \max_{\alpha >0 } \left\{-\frac{1}{n\alpha}+\sum_{i=1}^n\frac{1}{2d_i+\alpha} \right \}. \label{Delta-G}
\end{equation}
For arbitrary unweighted coupling graphs, the quantity \eqref{Delta-G} reduces to 
\begin{equation*}
		\Delta(\G) \hs = \hs -\frac{1}{2n} + \frac{n-1}{2n} \hs \sum_{i=1}^{n} \frac{1}{d_i},
	\end{equation*}
where the equality holds if  $\G$ is a complete graph  or complete bipartite graph.
\end{theorem}

For unweighted coupling graphs, the lower bound given by Theorem \ref{degree-2-thm} is tighter than the lower bound given by Theorem \ref{diam-thm}.  For $d$-regular  weighted coupling graphs, the lower bound is  
	\begin{equation}
	\Delta(\G) = \frac{(n-1)^2}{2nd}.
	\label{obtained}
	\end{equation}
This lower bound is tight for FOC networks with weighted coupling graphs, in the sense that the performance measure of a FOC network with  the weighted coupling graph $\K_n$ with identical edge weights $d/(n-1)$ meets the lower bound  (\ref{obtained}). 
%

%
%

	\begin{figure}[t]
		\centering
		\includegraphics[trim = 14mm 0mm 14mm 4mm, clip, width=0.7\textwidth]{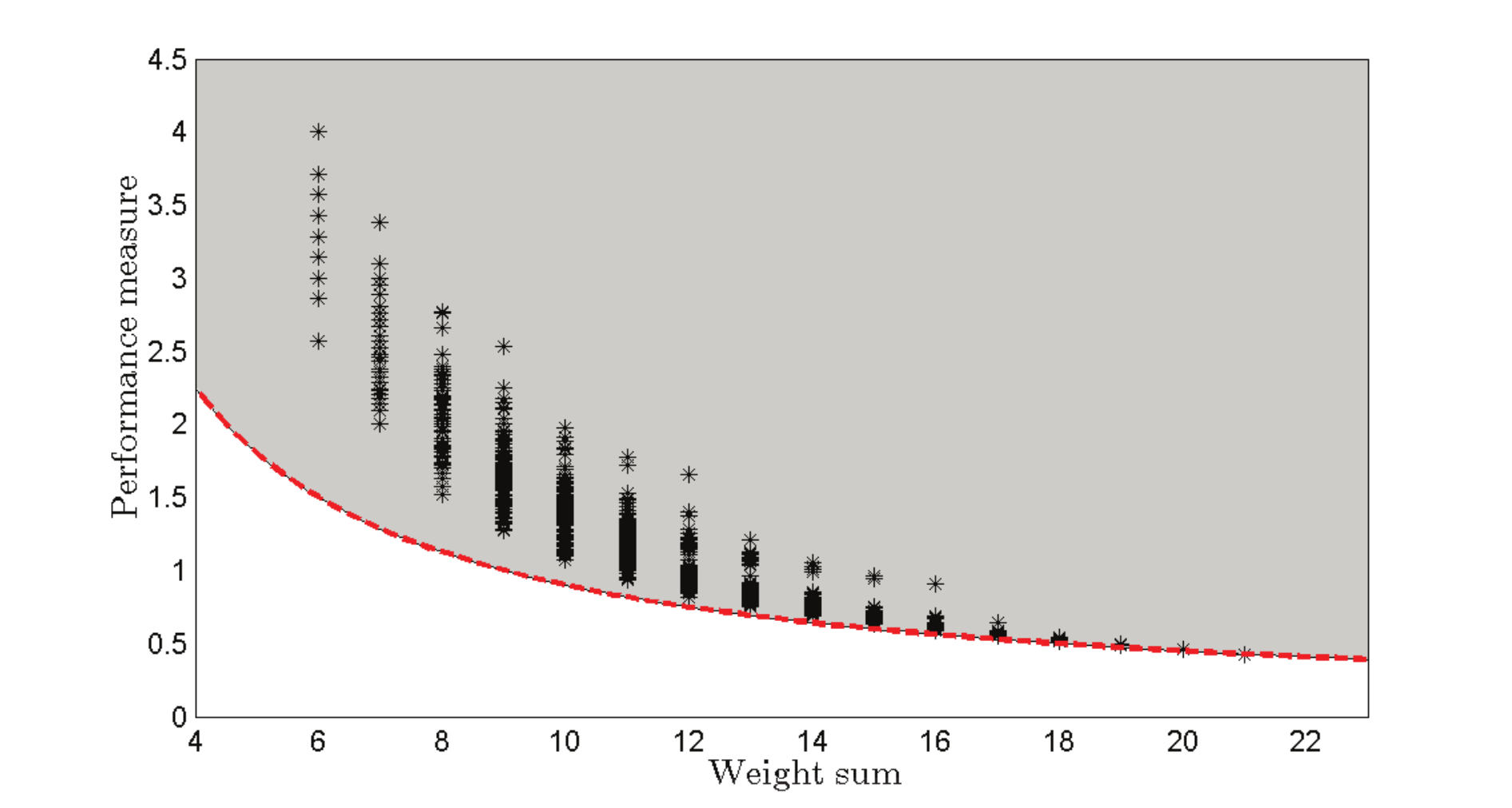}		
		\centering
  		\caption{{\small  This plot depicts the value of the performance measure for all FOC networks with coupling graphs in $\mathbb{G}_7$. The red dashed curve  portrays the lower bound in  \eqref{propo}.  }}
  		\label{plot1}
	\end{figure}

As we discussed earlier, our results for performance measure \eqref{rho-ss} can be utilized to analyze  performance measures \eqref{position-1} and \eqref{velocity-2} for SOC networks. In the following, we show that  there is an inherent lower bound on the best achievable performance measure \eqref{perf-meas-soc}.

\begin{theorem}\label{second-degree-thm}
For a given Type 2 SOC network \eqref{second-order}-\eqref{type-2} with an arbitrary weighted coupling graph in $\mathbb G_n$, the performance measure (\ref{perf-meas-soc}) is bounded from bellow by
	\begin{equation}
\RR (A^{(2)}_{\G}; M_n \oplus \mathbf 0) \hs \geq \hs		\frac{8m^4}{\beta \|L_{\G}\|_{F}^6}, \label{pos-lap-energy}
	\end{equation}
where $m$ is the number of edges in $\G$ and the Frobenius norm of the Laplacian matrix is defined by  
\begin{equation*}
\|L_{\G}\|_{F} := \left( \hs 2\sum_{e \in \E_{\G}}w_{\G}^2(e) \hs + \hs \sum_{i=1}^n d_i^{2} \hs \right)^{1/2}
\end{equation*}
and  $d_i$ for $i=1,\ldots, n$ are weighted node degrees. 
\end{theorem}

For $d$-regular unweighted graphs, the lower bound in \eqref{pos-lap-energy} reduces to 
	\begin{equation*}
\RR (A^{(2)}_{\G}; M_n \oplus \mathbf 0) \hs \geq \hs		\frac{n d}{2\beta(1+d)^3}, 
	\end{equation*}
and the equality holds if $\G=\mathcal K_n$.

\begin{remark}
In Theorem \ref{max-min-thm}, it is shown that the performance measure of a FOC network with an arbitrary unweighted coupling graph in $\mathbb{G}_n$  is always less than or equal to $(n^2-1)/12$. In the following, we show by means of three simple examples that the performance measure of a FOC network with a weighted  coupling graph can be made arbitrarily large. We consider a FOC network with three nodes and path coupling graph. The edge weights are given by $w(\{1,2\})=a$ and $w(\{2,3\})=1-a$, where $a > 0$. For different values of parameter $a$, the  total sum of edge weights is equal to $1$. However, we have $\RR (L_{\G};M_n)  \rightarrow \infty$ as $a \rightarrow 0$. 
Which implies that the performance measure cannot be uniformly bounded from above. Now for this graph, let us change the edge weights to  $w(\{1,2\})=a$ and $w(\{2,3\})=a^{-1}$.  According to \eqref{no-span-trees}, the total number of spanning trees of this graph is equal to $1$. It is straightforward to verify that $\RR (L_{\G};M_n) \rightarrow \infty$ as  $a \rightarrow 0$. In the third scenario, let us consider a cyclic graph with four nodes and edge weights $w(\{1,2\})=w(\{3,4\})=a$ and $w(\{2,3\})=w(\{1,4\})=1-a$. In this case, the weighted degree sequence is  $d_1=d_2=d_3=d_4=1$. A simple calculation shows that  $\RR (L_{\G};M_n) \rightarrow \infty$ as  $a \rightarrow 0$. These examples explain why the last column of Table \ref{table3} shows $\infty$ as the upper bound. 

\end{remark}

%

		\begin{figure}[t]
		\centering
		\includegraphics[trim = 14mm 0mm 14mm 4mm, clip, width=0.7\textwidth]{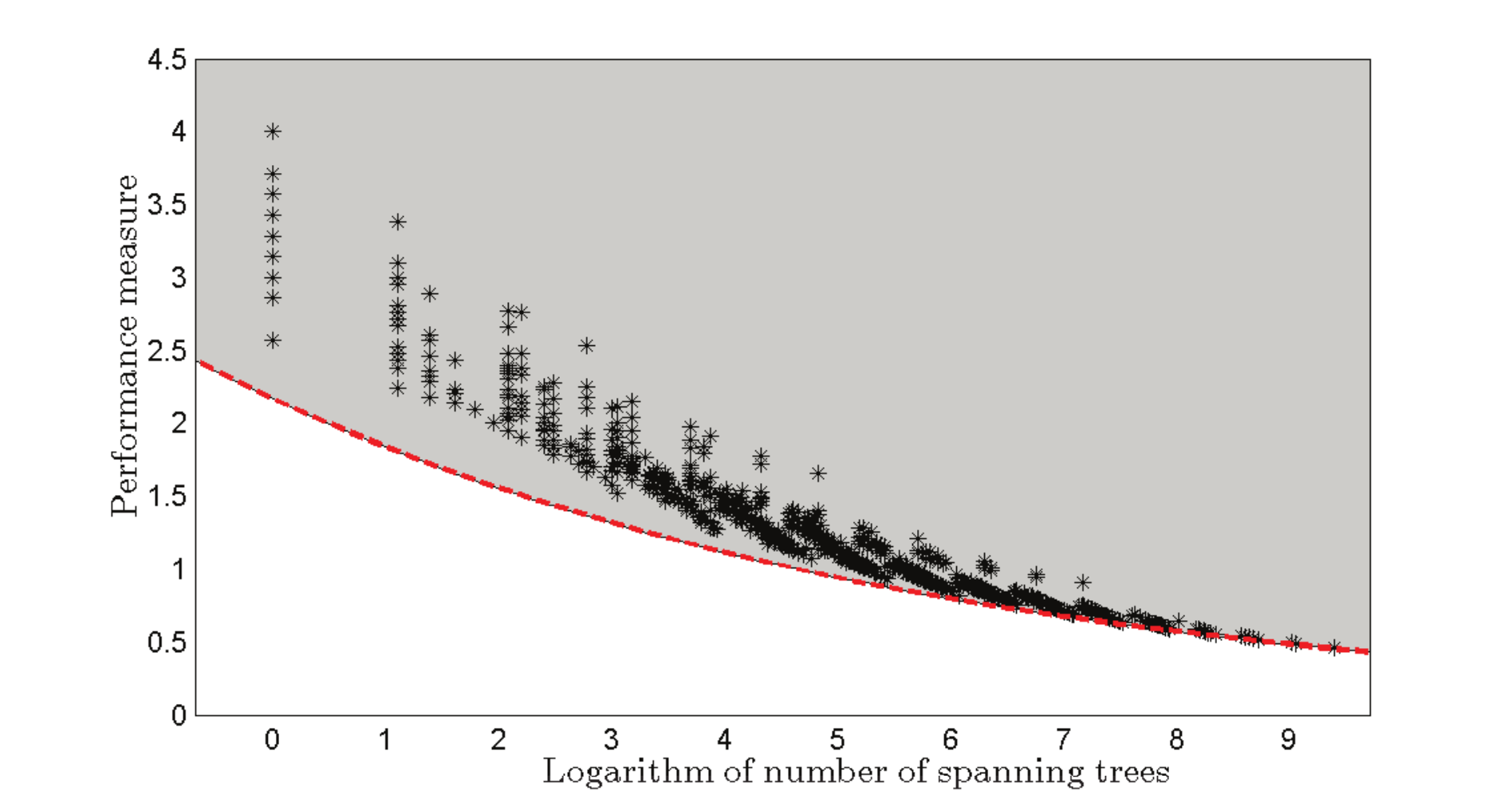}		
		\centering
  		\caption{{\small  This plot depicts the value of the performance measure for all FOC networks with coupling graphs in $\mathbb{G}_7$. The red dashed curve depicts the  lower bound in \eqref{ineq-spanning-tree}.}}
  		\label{plot2}
	\end{figure}

\subsection{Interpretation of Bounds as Fundamental Limits}\label{subsec-output-energy}

The value of the performance measure \eqref{perf-meas} for linear dynamical network \eqref{consensus-general} is equal to the average output energy of the autonomous linear dynamical network 
\begin{eqnarray}		
		\begin{cases}
		~\dot \psi~=~ -A \psi\\
		~y~=~ C_{\mathcal{Q}} \psi
		\end{cases}
		\label{consensus-general-1}
\end{eqnarray}
that is perturbed by an initial condition $\psi(0)$  drawn from an uncorrelated white stochastic process with the following characteristic   
\[ \mathbb E \left[ \psi(0) \psi(0)^{\text T} \right ] =  B B^{\text T}. \] 
In fact, it can be shown that  
\begin{equation}
\RR(A; L_{\Q}) \hs = \hs \mathbb E  \left[ \hs \int_0^{\infty} y(t;\psi(0))^{\text T}y(t;\psi(0))  \hs dt \hs \right],
\label{exp-y}
\end{equation}
where $y(t;\psi(0))$ is the output of the linear dynamical network \eqref{consensus-general-1} with respect to initial condition $\psi(0)$.  This  relationship enables us to equivalently interpret the performance measure \eqref{perf-meas} as the average energy needed to be consumed throughout  the network in order to steer the state of the randomly perturbed linear dynamical network to its equilibrium (i.e., consensus) state. Therefore, our theoretical bounds in Subsections  \ref{tailored-sub} and \ref{universal-sub} can be viewed as quantification of inherent fundamental limits on the minimum average energy required to be dissipated in the network in order to reach the consensus state again in steady state.

The use of term {\it fundamental} (or equivalently {\it hard}) limits for lower and upper bounds in our summary tables is appropriate and meaningful. The reason is that according to our results, the performance measure of a linear consensus network whose coupling graph has some known graph specification (e.g., number of nodes, number of spanning trees, total sum of edge weights, degree sequence, etc.)   cannot be better and worse than our theoretical lower bounds and upper bounds, respectively.

		\begin{figure}[t]
		\centering
		\includegraphics[trim = 14mm 0mm 14mm 4mm, clip, width=0.7\textwidth]{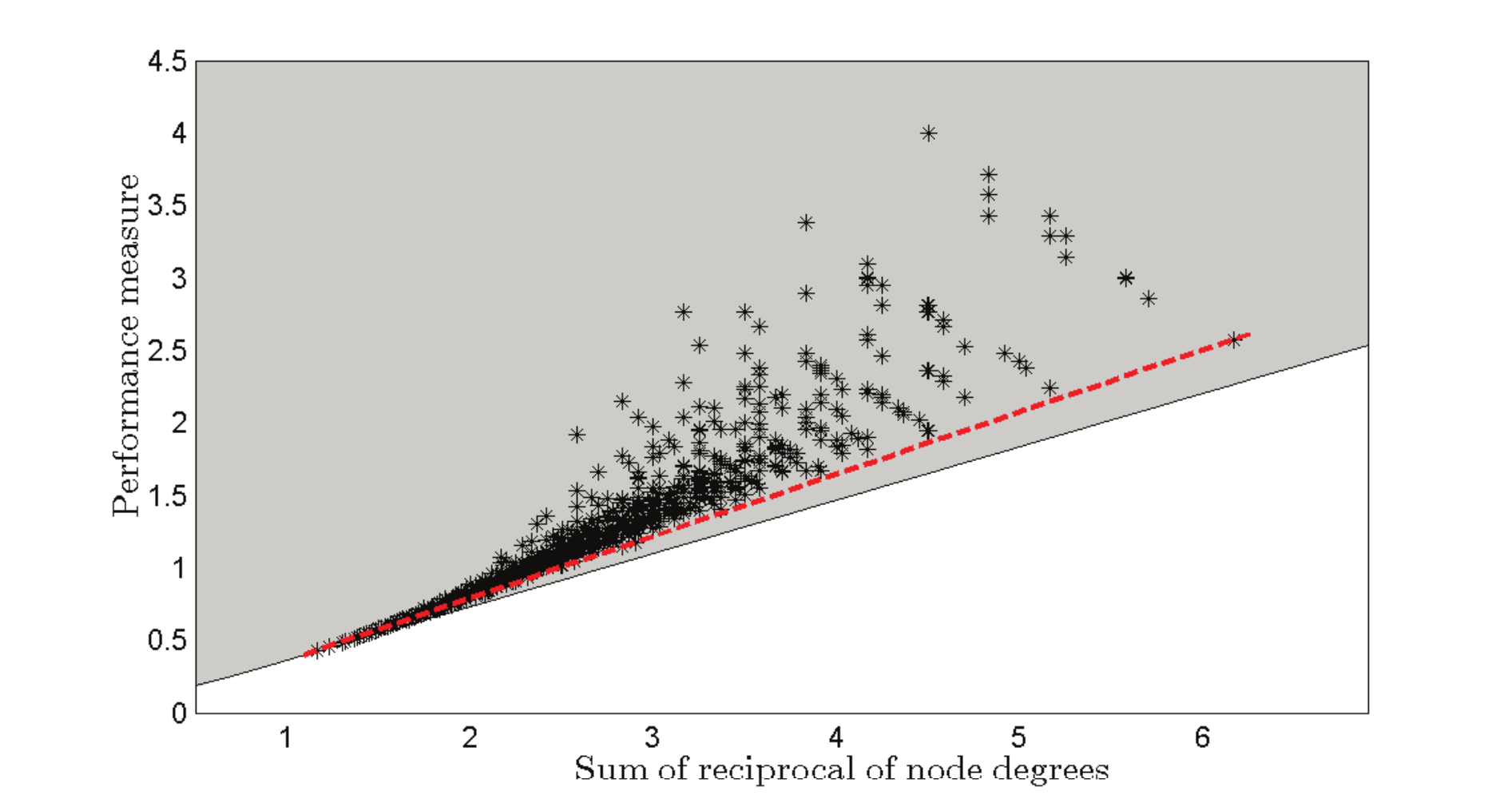}		
		\centering
  		\caption{{\small  This plot depicts the value of the performance measure for all FOC networks with coupling graphs in $\mathbb{G}_7$. The red dashed curve outlines the lower bound in \eqref{lower-delta} for unweighted graphs.}}
  		\label{plot3}
	\end{figure}

The philosophy behind our several results presented in Tables \ref{table1}, \ref{table2} and  \ref{table3} can be explained by portraying the value of performance measure for FOC and SOC networks versus various known graph specifications. In order to conceptualize the idea, we only focus on the class of FOC networks and three graph specifications in our analysis. Similar arguments can be extended for SOC networks and other graph specifications. Figures \ref{plot1}, \ref{plot2}, and \ref{plot3}  depict the value of the performance measures for FOC networks with coupling graphs in $\mathbb{G}_7$. In these figures, the points with star markers correspond to performance measures of all FOC networks with unweighted graphs in $\mathbb{G}_7$. The total number of such networks are $1,866,256$. In all three figures, the shaded grey area above  the red dashed curve corresponds to performance measures of FOC networks with weighted coupling graphs. In Figure \ref{plot1}, the performance measure \eqref{rho-ss} is drawn for different values of weight sum $W(\G)$. The lower bound in \eqref{propo} is highlighted by a red dashed curve and it draws a fundamental limit on the best achievable performance measures. One observes that the lower bound in \eqref{propo} is tight for a given value of weight sum. In fact, for a given $W(\G)$ there exists a weighted graph with total weight sum  $W(\G)$ whose performance measure reaches the exact value of the fundamental limit $\frac{(n-1)^2}{4W(\G)}$, where in this simulation $n=7$. However, this lower bound is conservative for unweighted graphs. For unweighted graphs, the weight sum is equal to the total number of edges in the coupling graph and it only assumes integer values. By exhausting all possible choices for unweighted graphs with identical number of edges in Figure \ref{plot1}, we show that there is a gap between the actual best achievable lower bound and our theoretical fundamental limit in \eqref{propo}. It can be perceived that this gap is  smaller for more dense coupling graphs. This observation suggests that our theoretical fundamental limit in \eqref{propo} is more conservative for sparse coupling graphs and less conservative for dense coupling graphs. Nevertheless, having more knowledge about graph specifications helps to close the gap. For example, the weight sums for FOC networks with tree and unicyclic coupling graphs are equal to $6$ and $7$, respectively. In these cases, the actual minimum and maximum achievable values of performance measure exactly matches with our theoretical fundamental limits in \eqref{uni-inq} and \eqref{tree-inq}. 

To conclude our discussion, one can also set out similar arguments for Figures \ref{plot2} and \ref{plot3} to infer that our theoretical fundamental limits in Subsection \ref{tailored-sub} are conservative for fairly sparse coupling graphs and much less conservative for dense coupling graphs. As we discussed in Subsection \ref{universal-sub}, one can exploit the structural properties of networks with sparse coupling graphs (e.g., trees and unicyclics) to quantify tight fundamental limits.

\section{Fundamental Tradeoffs Between Sparsity  and Performance Measure} 
\label{sec-V}

One of the design objectives for large-scale linear consensus networks is to optimize network coherence by designing a coupling graph that has the best possible sparsity and locality features. A fundamental property of performance measures \eqref{perf-meas-soc} and \eqref{rho-ss} is that they are monotonically decreasing functions of the coupling graphs in the cone of positive semidefinite matrices. This property implies that the value of the performance measure increases by sparsifying the underlying coupling graph, which is consistent with our results in Subsection \ref{universal-sub}. In this section, we quantify fundamental tradeoffs between the performance measure \eqref{rho-ss} and sparsity measures of FOC networks. The results of the following two theorems assert that the performance of a spanning subnetwork of a given FOC network never outperforms the performance of the parent network.

\begin{theorem} \label{order-thm}
Suppose that $\G \in \mathbb G_n$ is the coupling graph of a given FOC network. If $\mathcal P$ is a connected spanning subgraph of $\G$, then
	\begin{equation}
		\RR (L_{\G};M_n) \hs \leq  \hs \RR (L_{\mathcal P};M_n),
	\end{equation}
and the equality holds if and only if $\mathcal G = \mathcal P$.
\end{theorem}

\begin{theorem}\label{last-thm}
Suppose that $\G \in \mathbb G_n$ is the coupling graph of a given Type 2 SOC network. If $\mathcal P$ is a connected spanning subgraph of $\G$, then
\begin{eqnarray}
		\RR (A^{(2)}_{\G};  M_n \oplus \mathbf 0 ) & \leq & \RR (A^{(2)}_{\mathcal P};  M_n \oplus \mathbf 0), 					\label{subgraph-ineq} \\
		\RR (A^{(2)}_{\G}; \mathbf 0  \oplus M_n) & \leq & \RR (A^{(2)}_{\mathcal P}; \mathbf 0  \oplus M_n). \label{subgraph-ineq2} 
\end{eqnarray}
Moreover, the equalities hold if and only if $\mathcal G = \mathcal P$. 
\end{theorem}


The results of these two theorems implicitly assert that adding new edges to the underlying coupling graph of a consensus network may improve the global performance of the network. In the following, we identify several Heisenberg-like inequalities that quantify inherent fundamental tradeoffs betweens global performance and sparsity in FOC networks. The first sparsity measure that we consider is defined by 
	\begin{equation}
		\|\adj\|_{0} \hs := \hs \mathbf{card} \big\{ a_{ij} \neq 0~ | ~~\adj=[a_{ij}] \big\}, \label{sparsity-1}
	\end{equation}
where $\adj$ is the adjacency matrix of the coupling graph $\G$. For a given graph, the value of this sparsity measure is equal to twice the number of the edges. 
\begin{corollary}\label{1-coro}

For a given FOC network with an arbitrary  unweighted  coupling graph $\G \in \mathbb G_n$, there is a fundamental tradeoff between the performance measure \eqref{rho-ss} and the sparsity measure \eqref{sparsity-1} that is characterized in the multiplicative form by the following inequality
	\begin{eqnarray}
 		\RR (L_{\G};M_n) \hs \|\adj\|_{0} \hs \geq \hs \frac{(n-1)^2}{2}
		\label{tradeoff-1}
	\end{eqnarray}
and in the additive form by
	\begin{eqnarray}
		\frac{\RR (L_{\G};M_n)-\frac{1}{2}+\frac{1}{2n}}{\mathbf{diam}(\mathcal G)} ~+~ \frac{\|\adj\|_{0}}{4(n-1)} ~\leq~ \frac{n}{4}.
		\label{tradeoff-2}
	\end{eqnarray}
\end{corollary}

	Let us consider the class of  networks with identical number of nodes and compare several scenarios. The inequality \eqref{tradeoff-1} asserts that the best achievable levels of performance measure \eqref{rho-ss} for sparse FOC networks are comparably higher (worse) with respect to  less sparse FOC networks. For all FOC networks with identical diameters, inequality \eqref{tradeoff-2} implies that networks with more edges have smaller (better) levels of performance measures. Among all FOC networks with identical number of edges,  the ones with larger diameters have higher (worse) levels of performance measures.

\begin{corollary}\label{2-coro}
Let us consider the class of FOC networks with arbitrary unweighted  coupling graphs in  $\mathbb G_n$ and a given desired performance level $\rho_{ss}^*$. Then, the sparsity measure \eqref{sparsity-1} for this class of networks satisfies 
	\begin{equation}
		\frac{(n-1)^2}{2  \rho_{ss}^*} ~\leq~ \|\adj\|_{0} ~\leq~(n-1)\left[ n-4\left(\frac{\rho_{ss}^*-\frac{1}{2}+\frac{1}{2n}}{\mathbf{diam}(\mathcal G)} \right) \right]. \label{sparsity-bound}
	\end{equation}

\end{corollary}

The result of this corollary states that the graph diameter can be employed as a design parameter to achieve a desirable level of performance and sparsity. 

The second sparsity measure that we consider in this section is so called $\mathcal{S}_{0,1}$--measure and  defined by 
\begin{equation*}
		\|\adj\|_{\mathcal{S}_{0,1}} := \max \Big\{ \max_{1 \leq i \leq n} \| \adj(i,.) \|_{0}, \max_{1 \leq j \leq n} \| \adj(.,j) \|_{0}  \Big\},\end{equation*}
where $\adj(i,.)$ represents the $i$'th row and $\adj(.,j)$ the $j$'th column of adjacency matrix $\adj$.  The value of the $\mathcal{S}_{0,1}$--measure of a matrix is the maximum number of nonzero elements among all rows and columns of that matrix. We refer to \cite{motee-sun-IEEE-TAC-submitted-2014} for more details and discussions on this sparsity measure. The $\mathcal{S}_{0,1}$--measure of adjacency matrix of an unweighted graph is equal to the maximum node degree. The following result quantifies an inherent tradeoff between the performance measure and this sparsity measure. 

\begin{corollary}\label{3-coro}
For a given FOC network with an arbitrary  unweighted  coupling graph $\G \in \mathbb G_n$ with $n \geq 3$, there is a fundamental tradeoff between the performance measure \eqref{rho-ss} and the $\mathcal{S}_{0,1}$--measure that is characterized by 	\begin{equation}
		\left(\RR (L_{\G};M_n) + \frac{1}{2n}\right) \|\adj\|_{\mathcal{S}_{0,1}}  \hs \geq \hs \frac{n-1}{2}. \label{S-0-1-ineq}
	\end{equation}
\end{corollary}

The value of the $\mathcal{S}_{0,1}$--measure reveals some valuable information about sparsity as well as the spatial locality features
of a given adjacency matrix, while sparsity measure \eqref{sparsity-1} only provides information about sparsity. The inequality \eqref{S-0-1-ineq} asserts that the best achievable levels of performance measure \eqref{rho-ss}  decreases by improving local connectivity in the coupling graph of a FOC network.

The third sparsity measure of our interest for the class of FOC networks with unweighted coupling graphs is defined by 
	\begin{equation}
		\sigma(\G) := \max_{i, j \in \V_{\G} \atop i \neq j} \Big \{ d_i + d_j -  | \mathfrak{N}(i) \cap \mathfrak{N}(j) | \Big\},	
		\label{sss}
	\end{equation}
where $d_i$ is degree of node $i$ and $\mathfrak{N}(i)$ is the set of all nodes that are connected to node $i$ by an edge. 
The value of the sparsity measure $\sigma(\G)$ is equal to the maximum number of nodes that are connected to any pair of nodes among all pairs of nodes in the graph. It is easy to verify that  $\sigma(\G) \leq n$. The following result quantifies an inherent tradeoff between the performance measure and this sparsity measure. 

\begin{theorem}\label{13-thm}
For a given FOC network with an arbitrary  unweighted  coupling graph $\G \in \mathbb G_n$ with $n \geq 3$, there is a fundamental tradeoff between the performance measure \eqref{rho-ss} and sparsity measure \eqref{sss} that is quantified  by 	
	\begin{equation}
\RR (L_{\G};M_n) \hs \sigma(\G) \hs \geq \hs \frac{n-1}{2}. \label{eq1184}
	\end{equation}
Moreover, the equality holds if  $\G=\K_n$. 
\end{theorem}

To summarize our results in this section, we conclude that there are intrinsic fundamental tradeoffs between the two favorable design objectives in linear consensus networks: minimizing the performance measure and sparsifying the underlying coupling graph. 

{
\section{Performance Measure Interpretations for Two Real-World  Networks}\label{sec:applications}

In this section, we evaluate the performance measure for an interconnected power networks and a controlled group of vehicles in a formation. 

\subsection{Least Achievable Total Resistive Power Loss in Synchronous Power Networks}

In this part an example of Type 1 SOC networks is considered. Let us consider an interconnected network of synchronous generators with coupling graph $\mathcal{G}$ that consists of $n$ buses (nodes) and $m$ transmission lines (edges). A synchronous generator $G_i$ is associated to each node $i$ for $i=1,\ldots,n$ with inertia constant $M_i$, damping constant $\beta_i$, voltage magnitude $V_i$. It is assumed that a reduced order model of synchronous generator $G_{i}$ can be expressed using only two state variables: rotor angle $\theta_i$ and angular velocity $\omega_{i}$. Moreover, we assume that all damping constants are identical, i.e., $\beta=\beta_{1}=\ldots=\beta_{n}$.  For each edge $e \in  \E_{\mathcal G}$, we denote the admittance over $e$  by 
	\begin{equation}
		y_{e} = g_{e}-\mathbf j b_{e}, \label{admitt-edge}
	\end{equation}
where $g_{e}$ and $b_{e}$ are the conductance and susceptance of the corresponding transmission line, respectively, and $\mathbf{j}=\sqrt{-1}$. 
For each edge $e$, the ratio of its conductance to its susceptance is denoted by 
	\begin{equation}
		\alpha_e=\frac{g_{e}}{b_{e}}.
	\end{equation}

We define two graphs based on equation \eqref{admitt-edge}: conductance and susceptance graphs. The conductance graph is denoted and defined by $\mathcal{G}_g=(\V_{\G}, \E_{\G}, w_{\G_{g}})$ where $w_{\G_{g}}(e)=g_e$ for all $e \in \E_{\G}$. Similarly, the susceptance graph is denoted and defined by $\mathcal G_b=(\V_{\G}, \E_{\G}, w_{\G_{b}})$ where $w_{\G_{g}}(e)=b_e$ for all $e \in \E_{\G}$. In fact, the conductance and susceptance graphs are two identical copies of $\G$ but with different weight functions.

The governing nonlinear rotor dynamics of the interconnected network of synchronous generators (also known as swing equations) can be linearized  around the zero equilibrium operating point of the network in order to obtain
	\begin{eqnarray}
		&\left[ \begin{array}{ccc}
		\dot \theta \\
		\dot \omega \end{array} \right]=\left[ \begin{array}{ccc}
\mathbf 0& I  \\
		-L_{\G_b} & -\beta I  \end{array} \right] \left[ \begin{array}{ccc}
\theta \\
		\omega \end{array} \right]+\left[ \begin{array}{ccc}
\mathbf 0\\
		I \end{array} \right] \xi, \label{swing_1}
	\end{eqnarray}
where $\theta=\left[\begin{array}{ccc}\theta_{1} & \ldots & \theta_{n}\end{array}\right]^{\text T}$ and $\omega=\left[\begin{array}{ccc}\omega_{1} & \ldots & \omega_{n}\end{array}\right]^{\text T}$ are the state vectors of the entire network and $\xi$ is a zero-mean white noise process with identity covariance that models exogenous disturbances \cite{bamgay13acc,5530690}. 

The resistive power loss over each edge $e=\{i,j\}$ can be expressed as the following quantity 
	\begin{equation}
		P_{e}=g_{e}\hspace{0.03cm}|V_i-V_j|^2,
	\end{equation}
where $g_{e}$ is the conductance of edge $e$. Therefore, the total resistive power loss in the power network is given by  
	\begin{eqnarray}
		P_{\text{loss}}=\sum_{e=\{i,j\} \in \E_{\mathcal G}} P_{e}.
	\end{eqnarray}
If we consider the swing equations of the power network around its equilibrium point, we may apply the small angle approximation and replace the coupling terms $\sin(\theta_{i}-\theta_{j})$ by $\theta_{i}-\theta_{j}$ to obtain the following relationship
	\begin{eqnarray}
		\tilde P_{\text{loss}}=\sum_{e=\{i,j\} \in \E_{\mathcal G}} g_{e}\hspace{0.03cm}|\theta_i-\theta_j|^2.
		\label{p-loss}
	\end{eqnarray}
According to our definitions in Section \ref{sec:second-order}, the total resistive power loss $\tilde P_{\text{loss}}$ given by \eqref{p-loss} is equal to the performance measure of the linearized swing equation \eqref{swing_1} with respect to the angle output graph $\Q_{\theta}=\G_g$, where $\G_{b}$ is the corresponding conductance graph. Thus, according to Theorem \ref{2-thm}, we have
	\begin{equation}
		\RR \big( A_{\G_b}^{(1)};L_{\G_{g}} \oplus \mathbf 0 \big) \hs = \hs \frac{1}{2\beta}\mathbf{Tr}(L_{\G_b}^{\dag}L_{\G_g}).\label{H-2-loss}
	\end{equation}

In the following theorem, we show how the performance measure (\ref{H-2-loss}) can be expressed as a weighted mean of the ratios $\alpha_e$'s of all edges $e \in \mathcal E_\G$ in the network.


\begin{theorem}
\label{thm-power-loss}
The performance measure (\ref{H-2-loss}) of the linearized swing equations (\ref{swing_1}) with respect to the angle output graph $\G_{g}$ is given by 
	\begin{equation}
		\RR ( A_{\G_b}^{(1)};L_{\G_{g}}\oplus \mathbf 0 ) \hs = \hs \frac{\bar \alpha}{2\beta} (n-1),\label{total-loss}
	\end{equation}
and
	\begin{eqnarray}
		\bar \alpha=\frac{\sum_{e \in \E_{\mathcal G}} \nu_e \alpha_e}{\sum_{e \in \E_{\mathcal G}}\nu_e}=\frac{\sum_{e \in \E_{\mathcal G}}\nu_e \alpha_e}{n-1},
	\label{alpha-mean}
	\end{eqnarray}
in which $\nu_e=r^{(\G_b)}_e b_e$ and $r^{(\G_b)}_e$ and $b_{e}$ are the {\it effective reactance} and the susceptance of edge $e$, respectively. Furthermore, the expected  total resistive power loss is bounded by
	\begin{eqnarray}
		\frac{\alpha_{\min}}{2\beta}(n-1) \hs \leq \hs \RR ( A_{\G_b}^{(1)};L_{\G_{g}}\oplus \mathbf 0 ) \hs \leq \hs \frac{\alpha_{\max}}{2\beta} (n-1), 
		\label{loss-bounds}
	\end{eqnarray}
where  
	\begin{equation}
		\alpha_{\min} \hs = \hs \min_{{e \in \E_{\mathcal G}}} \alpha_e,~~~~~ \alpha_{\max} \hs = \hs \max_{e \in \E_{\mathcal G}} \alpha_e. 		
		\label{alpha-bar}
	\end{equation}
\end{theorem}
	
According to \eqref{total-loss}, the expected total resistive power loss depends on the specific structure of the coupling graph of the power network through $\bar {\alpha}$. However, the inequality \eqref{loss-bounds} shows that the lower and upper bounds of the expected total resistive power loss do not depend on the specific topology of the coupling graph of the network. For the special case when $\alpha_1 =\cdots=\alpha_m$, the result of Theorem \ref{thm-power-loss} asserts that the total power loss does not depend on the topology of the power grid (cf. \cite{bamgay13acc}). Under the assumption that all $\alpha_{e}$ are identical, the process of calculating the expected total resistive power loss benefits greatly from the symmetric structure of normal matrices \cite{Bamieh12}.

In the rest of this subsection for the case of nonidentical $\alpha_{e}$, we show that if the coupling graph has specific properties then the expected total resistive power loss is independent of the network topology but depends on the number of generators.

\begin{definition}
We say that graph $\mathcal G$ is an edge-transitive graph if there is an automorphism of $\mathcal G$ that maps $e_1$ to $e_2$ for all edges $e_1, e_2 \in \E_{\G}$. 
\end{definition}

Intuitively, in an edge-transitive graph all edges have identical local environments such that an edge cannot be distinguished from other edges based on its neighboring nodes and edges. Examples of edge-transitive graphs include biregular, star, cycle, $d$-dimensional torus $\mathbb Z^d_n$ and complete graphs \cite{Norman}.

\begin{theorem}
\label{edge-trans-thm}
Suppose that the coupling graph of the linearized power network (\ref{swing_1}) is edge-transitive and the internal susceptances  of all edges are identical. Then, the expected total resistive power loss is  given by
	\begin{equation}
		\RR ( A_{\G_b}^{(1)};L_{\G_{g}}\oplus \mathbf 0 ) \hs = \hs \frac{\sum_{e \in \E_{\mathcal G}}\alpha_{e}}{ 2\beta m}~ (n-1). 
		\label{price2}
		\end{equation}
	\end{theorem}

\begin{theorem}
\label{linear-loss-thm}
Suppose that the coupling graph of the linearized power network (\ref{swing_1}) is a tree $\G=\mathcal T $. Then, the expected total resistive power loss is  given by
	\begin{equation}
		\RR (A_{\mathcal T_b}^{(1)};L_{\mathcal T_{g}}\oplus \mathbf 0) \hs = \hs \frac{1}{ 2\beta}{\sum_{e \in \E_{\mathcal T}}\alpha_{e}}.
		\label{price1}
	\end{equation}
\end{theorem}
%

%
\subsection{Best Achievable Energy-Efficiency in Formation Control of Autonomous Vehicles}

The performance measures for a network of multiple autonomous vehicles have interesting output energy interpretations. Let us consider an abstract model of the formation control problem for a group of autonomous vehicles, which is given by a Type 2 SOC network in the form of \eqref{consensus-general-1} with state matrix \eqref{type-2}. Each vehicle has a position and a velocity variable and the state variable of the entire network is denoted by $\psi(t)=[\begin{array}{cc}x(t) & v(t)\end{array}]^{\text T}$ and is measured relative to a pre-specified desired trajectory $x_d(t)$ and velocity $v_d$. Without loss of generality, we assume that the position and velocity of each vehicle are scalar variables. The reason is that one can decouple  higher dimensional models into many Type 2 SOC models.  The overall objective is for the network to reach a desired formation pattern, where each autonomous vehicle travels at the constant desired velocity $v_d$ while preserving a pre-specified distance between itself and each of its neighbors. In this model, the state feedback controller uses both position and velocity measurements and $L_{\G}$ is in fact the corresponding feedback gain, which represents the coupling topology in the controller array, and constant $\beta$ is a design parameter  \cite{Bamieh11, Bamieh12}.  We assume that the initial condition of the overall network is drawn from a mutually uncorrelated white stochastic process that satisfies 
\[ \mathbb E \big[  \psi(0) \psi(0)^{\text T} \big] =  \left[\begin{array}{cc} \mathbf{0} & \mathbf{0} \\\mathbf{0} & I \end{array}\right],\] 
which implies that only the initial velocities are stochastically perturbed.  

As we discussed in Subsection \ref{subsec-output-energy}, the performance measure for linear consensus networks can be interpreted as the average output energy (\ref{exp-y}) required to be consumed in the network in order to steer the state of the randomly perturbed linear consensus network to its consensus state. The output energy with respect to complete position output graph is given by \eqref{perf-meas-soc} and it quantifies the average energy needed to be exhausted in the entire formation to follow the desired trajectory $x_d(t)$ in the steady state. The performance measure with respect to complete velocity output graph can be written as 
\begin{equation}
\RR (A_{\G}^{(2)};  \mathbf 0 \oplus M_n)~=~  \mathbb E \left[ \int_{0}^{\infty} \left( v^{\text T}(t)v(t)\hs - v^{\text T}_{\text{av}}(t)v_{\text{av}}(t) \right)\hs dt \hs \right], \label{eq1576}
\end{equation} 
where $v_{\text{av}}(t)=\frac{1}{n}J_n v(t)=\frac{1}{n}J_n v(0)$ is the (time-invariant) average velocity. The quantity in right hand side of \eqref{eq1576} represents the expected (extra) kinetic energy loss  in the network in order for all vehicles to reach the desired velocity $v_d$ in the steady state. The time-varying velocity of each vehicle leads to frequent accelerations and as result increases the vehicle's fuel consumption \cite{Mei2004}. Therefore, the performance measure \eqref{eq1576} closely depends on the total fuel consumption in the network. This energy interpretation implies that a fundamental limit on performance measure \eqref{eq1576} indicates the best achievable levels of energy-efficiency for a given network of autonomous vehicles in a formation.

\ignore{
	In the arguments above, we assumed that the velocity output graph is the centering graph. For a general velocity output graph with Laplacian matrix $L_{\Q_v}$, relation (\ref{eq1576}) is modified by:
	\begin{equation*}
		\RR (A_{\G}^{(2)};  \mathbf 0 \oplus L_{\Q_v}) \hs = \hs \mathbb{E} \left[ \int_0^{\infty} v^{\text T}(t){L_{\mathcal Q_v}}v(t) dt\right],
	\end{equation*}
where 
	\[v(t)^{\text T}{L_{\mathcal Q_v}}v(t) \hs = \hs \sum_{e=\{i,j\} \in \E_{\Q_v}} w_{\Q_v}(e)(v_i(t)-v_j(t))^2,\]
coincides with the energy of the flock according to \cite{4200853}. This value vanishes when all vehicles are moving with the same velocity. Therefore, the performance measure given by (\ref{velocity-2}) quantifies the steady state expected flock energy in presence of an exogenous white stochastic process.

}

\section{Discussion and Conclusion}

The primary focus of this paper is on a performance measure that is equal to the $\mathcal{H}_2$-norm (from exogenous disturbance input to an output) of linear consensus networks. The performance measures \eqref{perf-meas-soc} and \eqref{rho-ss} have several interesting functional properties. They are  convex  functions of Laplacian eigenvalues and monotonically decreasing with respect to adding new edges to the underlying coupling graph. The results of Section \ref{sec-V} highlights the importance of monotonicity property by quantifying inherent fundamental tradeoffs between sparsity  and the performance measures.  An interesting open problem to study whether these functional properties can lead us to categorize a larger class of admissible performance measures for linear consensus networks.

\newpage

\section{Appendix}
\section*{Definitions and Notations}
The following definitions and notations are used in our proofs in this appendix.

For a given Laplacian matrix $L_{\mathcal G}$, the corresponding resistance matrix $R_{\mathcal G}=[r_{ij}]$ is defined using the Moore-Penrose pseudo-inverse of $L_{\mathcal G}$ by setting $r_{ij}=l_{ii}^{\dag}+l_{jj}^{\dag}-l_{ji}^{\dag}-l_{ij}^{\dag}$, where $L^{\dag}_{\G}=[l^{\dag}_{ij}]$ and $r_{ij}$ is called the effective resistance between nodes $i$ and $j$. Furthermore, the total effective resistance $\mathbf r_{\textrm{total}}$ is defined as the sum of the effective resistances between all distinct pairs of nodes, i.e., 
	\begin{equation}
		\mathbf r_{\textrm{total}}~=~\frac{1}{2}\hspace{0.05cm} \mathbf{1}_{n}^{\text T} R_{\mathcal G} \mathbf{1}_{n}~=~\frac{1}{2}\sum_{i,j=1}^{n}r_{ij}.
	\end{equation}

We review some concepts from majorization theory. The following definition is from \cite{marshall11}.  

\begin{definition}
For every $x \in \R_+^n$, let us define $x^{\downarrow}$ to be a vector whose elements are a permuted version of elements of $x$ in descending order. We say that $x$  majorizes $y$, which is denoted by $x \unrhd y$, if and only if $\mathbf{1}^{\text T}x=\mathbf{1}^{\text T}y$ and 
	\begin{equation*}
		\sum_{i=1}^k x_i^{\downarrow} ~\geq~ \sum_{i=1}^k y_i^{\downarrow}, 
	\end{equation*}
for all $k=1,\ldots,n-1$.
\end{definition}

We should emphasize that majorization is not a partial ordering. This is because from relations $x \unrhd y$ and $y \unrhd x$ one can only conclude that the entries of these two vectors are equal, but not necessarily in the same order. Therefore, relations $x \unrhd y$ and $y \unrhd x$ do not imply $x=y$.  

%
\begin{definition}
The real-valued function $F: \R_+^n \rightarrow \R$ is called Schur--convex if $F(x) \geq F(y)$ for every two vectors  $x$ and $y$ with property $x \unrhd y$. Similarly, a function $F$ is Schur--concave if $-F$ is Schur--convex.
\end{definition}

\section*{Proof of Theorem \ref{1-thm}}
\label{proof-1-thm}
Here we present two approaches to calculate the value of the performance measure:

{\noindent \bf First proof:}
 Let us define the disagreement vector by (cf. \cite{olfati})
	\begin{equation}
		x_d(t) ~:=~ M_n x(t)~=~x(t) - \frac{1}{n}J_nx(t).
		\label{dis-vector-1}
	\end{equation}
By multiplying a vector by the centering matrix, we actually subtract the mean of all the entries of the vector from each entries. The dynamics of $\NN(L_{\G};L_{\Q})$ with respect to the new state transformation \eqref{dis-vector-1} is so called disagreement form of the network, which is given by
	\begin{equation}
		\NN_{d}(L_{\G_{d}};L_{\Q}):\begin{cases}
		\dot x_d=-L_{\G_d} x_d + M_n \xi \\
     		y~=~C_{\Q} \hspace{0.05cm}x_d
		\end{cases}
		\nonumber
	\end{equation}
in which the new system matrix is $-L_{\G_d}=-(L_{\mathcal G}+\frac{1}{n}J_n)$ and it is stable. One can easily verify that the transfer functions from $\xi$ to $y$  in both networks $\NN(L_{\G};L_{\Q})$ and $\NN_{d}(L_{\G_{d}};L_{\Q})$ are identical. Therefore, the $\mathcal H_2$--norm of the system from $\xi$ to $y$ in both representations are well-defined and equivalent. Therefore, we consider the integral form of the output of network $\NN_d(L_{\G_d};L_{\Q})$ as follows
\begin{equation}
	y(t)~=~C_{\Q} \int_{0}^{t} e^{-L_{\G_d}(t-\tau)} M_n \xi(\tau) \hs d\tau. \label{output} 
\end{equation}
By substituting $y(t)$ from \eqref{output} in \eqref{perf-meas}, calculating the expected  value, and finally taking the limit, the value of the performance measure can be calculated using the trace formula $\Tr(P_cL_{\Q})$, where matrix  $P_c$ is the controllability Gramian of the disagreement network $\NN_d(L_{\G_d};L_{\Q})$ and it is the solution of the Lyapunov equation
\begin{equation*}
	L_{\G_d} P_c+P_c L_{\G_d} - M_n=0.
\end{equation*}  
Note that in this case $-L_{\G_d}$ is stable, therefore this Lyapunov equation has a unique solution \cite[Th. 7.11]{Rugh}.
Using $L^{\dag}_{\G}L_{\G_d}=L_{\G_d}L^{\dag}_{\G}=M_n$, we get $P_c=\frac{1}{2} L^{\dag}_{\G}$. Therefore, we get our desired result.

{\noindent \bf Second proof:} According to \cite{Doyle89}, we need to calculate 
	\begin{equation}
		\RR(L_{\G}; L_{\Q})~=~ \frac{1}{2\pi} \int_{-\infty}^{\infty} \Tr ( G^{*}(j\omega)G(j\omega)) 
		\hspace{0.05cm} d\omega,
		\label{h2norm}
	\end{equation}
where $G(s)$ is the transfer function of the FOC network (\ref{first-order-1}) from $\xi(t)$ to $y(t)$
	\begin{equation*}
		G(s)~=~C_{\Q}(sI+L_{\G})^{-1}.
	\end{equation*}
Then we rewrite the integrand of (\ref{h2norm}) in the following form
	\begin{equation}
		\Tr ( G(j\omega)G^{*}(j\omega))~=~ \Tr (C_{\Q}M_n( \omega^2 I+L_{\G}^2)^{-1}M_nC_{\Q}^{\text T}),
		\label{integrant} 
	\end{equation}
where in (\ref{integrant}) we use the fact that $C_{\Q}=C_{\Q}M_n$. From (\ref{h2norm}) and (\ref{integrant}), it follows that
	\begin{eqnarray}
		\RR(L_{\G}; L_{\Q}) \nonumber &=& \frac{1}{2\pi} \int_{-\infty}^{\infty} \Tr \left (L_{\Q}M_n( \omega^2 I+L_{\G}^2)^{-1}M_n\right )\hspace{0.05cm} d\omega \nonumber \\
		~~&=& \frac{1}{2\pi} \Tr \left( L_{\Q} \int_{-\infty}^{\infty} M_n( \omega^2 I+L_{\G}^2)^{-1}M_n \hspace{0.05cm} d\omega \right).
		\label{integrant2}
	\end{eqnarray}
We now consider the eigenvalue decomposition of the Laplacian matrix which is given by 
	\begin{equation*}
		L_{\G}~=~ U\Lambda  U^{\text T},
	\end{equation*}
where $\Lambda = \mathbf{diag}(\lambda_1,\ldots,\lambda_{n})$ and $ U = [\mathbf{u}_1, \mathbf{u}_2, \cdots, \mathbf{u}_n]$ is the corresponding orthonormal matrix of eigenvectors. Consequently, we can rewrite $M_n$ and $L_{\G}$ as follows
	\begin{equation*}
		M_n~=~U \mathbf{diag} \left ( 0,1,\cdots, 1\right )U^{\text T},
	\end{equation*}
and
	\begin{equation*}
		L_{\G}~=~U \mathbf{diag} \left ( 0,\lambda_2,\cdots, \lambda_n\right )U^{\text T}.
	\end{equation*}
Note that graph $\G$ is connected, therefore the corresponding Laplacian matrix has only one zero eigenvalue $\lambda_1=0$ with corresponding eigenvector $\mathbf 1$.   
The integrand of (\ref{integrant2}) can be rewritten as 
	\begin{eqnarray}
		&&\hspace{-1cm}M_n( \omega^2 I+L_{\G}^2 )^{-1}M_n~=~ U\mathbf{diag} \left (0,\frac{1}{\omega^2+\lambda_2^2},\cdots, \frac{1}{\omega^2+\lambda_n^2}\right )U^{\text T}.
		\label{decom2}
	\end{eqnarray}
Finally, from (\ref{decom2}) and (\ref{integrant2}), we get 
	\begin{eqnarray}
		\RR(L_{\G}; L_{\Q})
		&=& \frac{1}{2\pi} \Tr ( L_{\Q} U\mathbf{diag} \left (0,{\pi}{\lambda_2^{-1}},\cdots, {\pi}{\lambda_n^{-1}}\right )U^{\text T} )\nonumber \\
		&=& \frac{1}{2} \Tr ( L_{\Q} L_{\G}^{\dag} ).
	\end{eqnarray}

\section*{Proof of Theorem \ref{2-thm}}
\label{proof-2-thm}
Similar to the proof of Theorem \ref{1-thm}, one can define the disagreement position and velocity vectors by 
	\begin{equation*}
		x_d(t) ~:=~ M_n x(t)~~~\text{and}~~~v_d(t)~:=~M_n v(t).
	\end{equation*}
and then follow the same approach presented in the proof of Theorem \ref{1-thm}. Calculating the controllability Gramian matrix of the new disagreement SOC network is straightforward.
%

\section*{Proof of Theorem \ref{max-min-thm}}
\label{proof-max-min-thm}
Theorem \ref{order-thm} implies that for any graph $\G$ with $n$ nodes, we have
	\[ \RR (L_{\G};M_n) ~ \geq ~ \mathbf  \RR (L_{\K_n};M_n),\]
because graph $\G$ is always a subgraph of $\K_n$; then an explicit computation shows that $ \RR (L_{\K_n};M_n)=(n-1)/(2n)$.
On the other hand, $ \RR $ reaches its maximal value when the coupling graph is a tree. We refer to Theorem \ref{tree-thm} for more details and a proof.

\section*{Proof of Theorem \ref{tree-thm}}
\label{proof-tree-thm}
Consider the  characteristic polynomial of the Laplacian matrix of the coupling graph $\mathcal T$
	\begin{equation}
		\Phi_{\mathcal T}(\lambda)~=~\sum_{k=0}^n \hspace{0.05cm}(-1)^{n-k}\hspace{0.05cm}c_k(\mathcal T)\hspace{0.05cm}\lambda^k.
		\label{eq-pol-tree}
	\end{equation}
From (\ref{eq-h}) and Vieta's formulas for (\ref{eq-pol-tree}), it follows that
	\begin{equation}
		\RR \left(L_{\mathcal T}; M_n \right)~=~\frac{c_2(\mathcal T)}{2 c_1(\mathcal T)}.
		\label{eq-h-co-tree}
	\end{equation}
We also know that $c_1(\mathcal T)=\prod_{i=2}^n \lambda_i$; and this quantity is equal to $n$ for tree graphs. Therefore, one can rewrite (\ref{eq-h-co-tree}) as follows
	\begin{equation}
		\RR \left(L_{\mathcal T}; M_n \right)~=~\frac{c_2(\mathcal T)}{2 n}.
		\label{eq--tree}
	\end{equation}
One of the invariant characteristics of a graph is its Wiener number which is denoted by $\mathbf W(\mathcal{T})$\cite{Gutman2003}. This quantity is equal to the sum of distances between all pairs of nodes of $\mathcal T$. It is well known that the second coefficient of the Laplacian characteristic polynomial of a tree coincides with the Wiener number, \ie 
	\[c_2({\mathcal T})~=~\mathbf W(\mathcal T).\] 
According to this fact and \eqref{eq--tree}, it follows that
	\begin{equation}
		\RR ({\mathcal T})~=~\frac{\mathbf W(\mathcal T)}{2n}.
		\label{h-wiener}
	\end{equation}
Based on reference \cite{tree} if $\mathcal T$ is a tree with $n$ nodes that is neither $\mathcal P_n$ nor $\mathcal S_n$, then
	\begin{eqnarray}
		\mathbf W(\mathcal S_n)~<~\mathbf W(\mathcal T)~<~\mathbf W(\mathcal P_n).
		\label{w}
	\end{eqnarray}
Furthermore, it is shown that (cf. \cite{tree})
	\begin{eqnarray}
		\mathbf W(\mathcal P_n)=\binom{n+1}{3},~~\text{and}~~\mathbf W(\mathcal S_n)=(n-1)^2.
		\label{w-com}
	\end{eqnarray}
From  (\ref{h-wiener}), (\ref{w}) and (\ref{w-com}), we have  
	\[\frac{(n-1)^2}{2n} ~<~ \RR (L_{\mathcal T}; M_n) ~<~ \frac{n^2-1}{12}.\]
On the other hand, it follows from (\ref{w-com}) and (\ref{h-wiener}) that
	\begin{equation*}
		\small{\RR (L_{\mathcal P_n};M_n)=\frac{n^2-1}{12},~~\text{and}~~\RR (L_{\mathcal S_n};M_n)=\frac{(n-1)^2}{2n}}.
	\end{equation*}
Therefore,  the lower bound in (\ref{tree-inq}) is achieved if and only if  $\mathcal{T}=\mathcal S_n$, and the upper bound is achieved if and only if  $\mathcal{T}=\mathcal P_n$.
\section*{Proof of Theorem \ref{unicyclic-thm}}
\label{proof-unicyclic-thm}

The proof is a direct consequence of \cite[Th. 4.4]{unicyclic} and \eqref{eq-r}.

\section*{Proof of Theorem \ref{bipartite-thm}}
\label{proof-bipartite-th} 
	According to Theorem \ref{max-min-thm}, a path graph $\mathcal P_n$ has the maximal level of performance measure among all graphs with $n$ nodes. Moreover,  $\mathcal P_n$ is in fact a bipartite graph. Therefore, we get
	\begin{equation}
		\RR (L_{\mathcal G}; M_n) ~\leq~ \frac{n^2-1}{12}.
		\nonumber
	\end{equation}
The best achievable lower bound can be obtained  from (\ref{eq-r}) and the result of \cite[Th. 3.1]{Yang}.

\section*{Proof of Theorem \ref{general-thm-type-2}}
\label{proof-second-degree-thm-3}
Consider the  characteristic polynomial of the Laplacian matrix of the coupling graph $\G$
	\begin{equation}
		\Phi_{\G}(\lambda)~=~\sum_{k=0}^n \hspace{0.05cm}(-1)^{n-k}\hspace{0.05cm}c_k(\G)\hspace{0.05cm}\lambda^k.
		\label{eq-pol}
	\end{equation}
From (\ref{perf-meas-soc}) and Vieta's formulas for (\ref{eq-pol}), it follows that
	\begin{equation}
		\RR (A^{(2)}_{\G}; M_n \oplus \mathbf 0) ~=~\frac{1}{2\beta}\left [\left(\frac{c_2(\G)}{c_1(\G)}\right)^2-\frac{c_3(\G)}{c_1(\G)}\right].
		\label{eq-h-co}
	\end{equation}
We also know that $c_1(\G)=n \mathfrak T(\G)$, therefore we can rewrite (\ref{eq-h-co}) as follows
	\begin{equation}
		\RR (A^{(2)}_{\G}; M_n \oplus \mathbf 0) ~=~\frac{c_2^2(\G)-2nc_3(\G)}{2\beta n^2}.
		\label{eq-h-co}
	\end{equation}
Based on Theorem \ref{last-thm}, for finding the upper bound on the performance measure among all graphs in $\G_n$, we can only focus on tree graphs.
According to \cite[Th.1]{Gutman2003}, for a given tree with $n\geq 5$ nodes and different from path and star, we have 
\[ (n-1)^2 = c_2(\mathcal S_n) ~<~ c_2(\mathcal T) ~< ~c_2 (\mathcal P_n) = \binom{n+1}{3},\]
and
\[ \frac{(n-1)(n-2)^2}{2} = c_3(\mathcal S_n) ~<~ c_3(\mathcal T) ~< ~c_3 (\mathcal P_n) = \binom{n+2}{5}.\]
Using these equations and \eqref{eq-h-co}, we get 
	\begin{eqnarray*}
		\RR (A^{(2)}_{\G}; M_n \oplus \mathbf 0) &<& \frac{c_2^2(\mathcal P_n)-2nc_3(\mathcal S_n)}{2\beta n^2} \nonumber \\
		&<&\frac{(n^2-1)^2}{72\beta}-\frac{(n-1)(n-2)^2}{2n\beta}.
	\end{eqnarray*}
On the other hand, according to Theorem \ref{last-thm}, the desired lower bound is achieved for a complete graph $\G=\K_n$.
\section*{Proof of Theorem \ref{tree-thm-type-2}}
\label{proof-second-degree-thm-1}

Similar to the proof of Theorem \ref{general-thm-type-2}, we can get the desired upper bound. On the other hand, for an unweighted tree graph $\mathcal T \neq \mathcal S_n$ with $n$ nodes, we have 
	\begin{equation}
		1~=~\lambda_2(\mathcal S_n) ~>~ \lambda_2(\mathcal T).
		\label{eq2282}
	\end{equation}
We also know that $\sum_{i=2}^n \lambda_i(\mathcal T)=2n-2$. Therefore, we get
	\[ \begin{bmatrix}1& \frac{2n-3}{n-1}& \ldots & \frac{2n-3}{n-1} \end{bmatrix}^{\text T} ~\unlhd ~ \begin{bmatrix} \lambda_2(\mathcal T)& \cdots &\lambda_n(\mathcal T)\end{bmatrix}^{\text T}.\] 
Moreover, it can be easily shown that $\RR \left(A^{(2)}_\G; M_n \right)$ is a Schur--convex function respect to $\begin{bmatrix}\lambda_2(\mathcal T)& \cdots & \lambda_n(\mathcal T)\end{bmatrix}^{\text T} \in \R^{n-1}_{++}$, where $\R_{++}$ denotes the set of all positive real numbers. Therefore, according to the definition of Schur--convex functions, one can obtain the desired lower bound. 

\section*{Proof of Theorem \ref{diam-thm}}
\label{proof-diam-thm}

For the lower bound, we apply the inequality of arithmetic and harmonic means and (\ref{eq-h})
	\begin{eqnarray*}
		\RR \left(L_{\G};M_n\right)~=~\frac{1}{2}\sum_{i=2}^n\lambda_i^{-1}~\leq~\frac{n-1}{2 \sum_{i=2}^n \lambda_i} ~=~ \frac{n-1}{4 W(\G)}.
	\end{eqnarray*}
%
%
On the other hand, according to (\ref{eq-r}) for the upper bound we get
    \begin{equation}
        \RR \left(L_{\G};M_n \right)~=~\frac{1}{2n}\sum_{i\neq j} r_{\{i,j\}}~=~\frac{1}{2n}\left (\sum_{e \in \E_{\G}}r_e+\sum_{e \notin \E_{\G}} r_e\right ).
        \label{upper}
    \end{equation}
Moreover, based on \cite[Lemma 2]{Rapat} for unweighted graph we have
    \begin{equation}
        \sum_{e \in \E_{\G}} r_e~=~ n-1.
        \label{n-1}
    \end{equation}
From (\ref{upper}) and (\ref{n-1}), it follows that
    \begin{equation}
        \RR \left(L_{\G};M_n \right)~=~\frac{n-1}{2n}~+~\frac{1}{2n}\sum_{e \notin \E_{\G}} r_e.
        \label{up}
    \end{equation}
We note that the distance between two nodes of graph $\G$ is less than or equal to $\mathbf{diam}(\G)$, therefore $r_{\{i,j\}} \leq \mathbf{diam}(\G)$, using this fact and (\ref{up}), we get the desired upper bound
    \begin{equation*}
        \RR \left(L_{\G};M_n \right)~\leq~\frac{1}{2n}\left(n-1~+~\left[\binom{n}{2}-m\right]\mathbf{diam}(\G)\right).
    \end{equation*}
%


\section*{Proof of Theorem \ref{trace-thm}}
\label{proof-trace-thm}

It can be shown that $\RR (L_{\G};M_n)$ is  a Schur--convex function with respect to $[\lambda_2,\ldots, \lambda_n]^{\text T} \in \R^{n-1}_{++}$ where $\lambda_{i}$ for $i=2,\ldots,n$ are eigenvalues of $L_{\G}$. On the other hand, we have 
	\[\frac{\mathbf{Tr}(L_{\G})}{n-1} \hs \mathbf 1_{n-1}^{\text T} ~\unlhd~ \begin{bmatrix} \lambda_2& \cdots &\lambda_n \end{bmatrix}^{\text T}.\] 
Therefore, according to the definition of Schur--convex functions, we can conclude inequality (\ref{propo}).

\section*{Proof of Theorem \ref{spanning-thm}}
\label{proof-spanning-thm}

By applying the inequality of arithmetic and geometric means to (\ref{eq-h}), we get
	\begin{equation}
		{ \RR (L_{\G};M_n)}~=~\frac{1}{2}\sum_{i=2}^n \lambda_i^{-1} ~\geq ~ \frac{n-1}{2}\sqrt[n-1]{\prod_{i=2}^n \lambda_i^{-1}}.
		\label{in-3}
	\end{equation}
Using Kirchhoff's matrix tree theorem the number of spanning trees of graph can be expressed as follows 
	\begin{equation}
		\mathfrak{T}(\mathcal G)~=~\frac{1}{n}\prod_{i=2}^n \lambda_i.
		\label{in-4}
	\end{equation}
Using (\ref{in-3}) and (\ref{in-4}), we get the desired lower bound.

\section*{Proof of Theorem \ref{cut-thm}}

It is shown that the performance measure of FOC network (\ref{first-order-1}) can be calculated by \eqref{eq-r}.
Moreover, in reference \cite{Deng}, it is shown that the $\mathbf r_{\textrm {total}}$ can be bounded from below as 
	\begin{equation*}
		\mathbf r_{\textrm {total}} ~\geq~ n\left (\kappa(\G)+1\right)+1-\frac{2n}{n-\kappa(\G)},
		\label{r-tot}
	\end{equation*}
for all connected graphs with $n$ nodes and $\kappa(\G)$ cut edges. The lower bound can be achieved if and only if $\mathcal G=\mathcal S(\mathcal K_{n-\kappa({\G})}; \mathcal K_1, \cdots , \mathcal K_1)$. 


\section*{Proof of Theorem \ref{degree-2-thm}}
\label{proof-degree-2-thm}
We consider two cases:\\
{\bf Weighted graph:} Assume that $\tilde L_{\G}= L_{\G}+{\alpha} J_n$ and $\alpha>0$, note that the eigenvalues of $\tilde L_{\G}$ are $n\alpha, \lambda_2, \cdots, \lambda_n$, where $\lambda_i$'s are eigenvalues of $L_{\G}$. Based on Schur--Horn theorem the diagonal elements of $\tilde L_{\G}$ are majorized by its eigenvalues, therefore we have
	\begin{equation}
		\sum_{i=1}^n\frac{1}{d_i+\alpha}~ \leq~ \frac{1}{n\alpha}+\sum_{i=2}^n \lambda_i^{-1}.
		\label{Schur}
	\end{equation}	
From the definition of $\RR(L_{\G};M_n)$ and (\ref{Schur}), it follows that
	\begin{equation}
		\frac{-1}{n\alpha}+\sum_{i=1}^n\frac{1}{2d_i+\alpha}~ \leq ~\RR \left(L_{\G};M_n\right).
		\label{Schur-2}
	\end{equation}	
{\bf Unweighted graph:} Using the same idea in the proof of Theorem \ref{diam-thm}, we can rewrite the performance measure of a FOC network (\ref{first-order-1}), as follows
    \begin{equation*}
        \RR (L_{\G};M_n)~=~\frac{n-1}{2n}~+~\frac{1}{2n}\sum_{e \notin \E_{\G}} r_e.
        \label{up22}
    \end{equation*}
Note that $r_{\{i,j\}} \geq \frac{1}{d_i}+\frac{1}{d_j}$, this implies
    \begin{eqnarray}
        \RR (L_{\G};M_n)&\geq&\frac{n-1}{2n}~+~\frac{1}{2n}\sum_{\{i,j\} \notin \E_{\G}} (\frac{1}{d_i}+\frac{1}{d_j})\nonumber \\
        &=&\frac{n-1}{2n}~+~\frac{1}{2n}\sum_{i=1}^n\frac{n-1-d_i}{d_i}\nonumber \\
        &=&\frac{-1}{2n}~+~\frac{n-1}{2n}\sum_{i=1}^n\frac{1}{d_i}.\nonumber 
        \label{up22-2}
    \end{eqnarray}
This completes the proof. The interested reader is referred to \cite{Zhou2} for more details and similar arguments.

\section*{Proof of Theorem \ref{second-degree-thm}}
\label{proof-second-degree-thm-2}

From the H\"older's inequality, it follows that
	\begin{eqnarray}
		\sum_{i=2}^n\lambda_i&=&\sum_{i=2}^n\left( \lambda_i^{-\frac{1}{2}}\right)\left (\lambda_i^{\frac{3}{2}}\right ) \nonumber\\
		&\leq&\left (\sum_{i=2}^n\left (\lambda_i^{-\frac{1}{2}}\right)^4\right)^{\frac{1}{4}}\left (\sum_{i=2}^n \left (\lambda_i^{\frac{3}{2}}\right )^{\frac{4}{3}} \right)^{\frac{3}{4}} \nonumber\\
		&=&\left (\sum_{i=2}^n\lambda_i^{-2}\right)^{\frac{1}{4}}\left (\sum_{i=2}^n \lambda_i^{2} \right)^{\frac{3}{4}}.
		\label{H_(2)-ineq}
	\end{eqnarray}
The inequality (\ref{H_(2)-ineq}) can be rewritten in the following form 
	\begin{eqnarray}
		\frac{\sum_{i=2}^n\lambda_i}{\left (\sum_{i=2}^n \lambda_i^{2} \right)^{\frac{3}{4}}} ~\leq~ \left (\sum_{i=2}^n \lambda_i^{-2}\right)^{\frac{1}{4}}.
		\label{eq2352}
	\end{eqnarray}
By combining (\ref{perf-meas-soc}) and (\ref{eq2352}) and using the facts that $\sum_{i=2}^n \lambda_i=2m$ and $\|L_{\mathcal G}\|_{F}^2=\sum_{i=2}^n \lambda_i^2$, we have
	\begin{equation*}
		\frac{2m}{\|L_{\mathcal G}\|_{F}^{1.5}} ~\leq~ \left (2 \beta \RR \left(A^{(2)}_{\G};  M_n \oplus \mathbf 0  \right)\right)^{\frac{1}{4}}.
	\end{equation*}
Thus, one can conclude that \eqref{pos-lap-energy} holds.

\section*{Proof of Theorem \ref{order-thm}}
\label{proof-order-thm}

For every $x \in \R^n$, we have
	\begin{eqnarray}
		x^{\text T}L_{\mathcal G}x&=&\sum_{e=\{i,j\} \in \E_{\mathcal G}} w(e) \left (x_i-x_j \right )^2 \nonumber \\ 
		&\geq& \sum_{e=\{i,j\} \in \E_{\mathcal P}} w(e) \left (x_i-x_j \right )^2 \nonumber\\
		&= & x^{\text T}L_{\mathcal P}x.
		\label{ineq-tr1}
	\end{eqnarray} 
This inequality implies that $L_{\mathcal P} ~ \leq ~ L_{\mathcal G}$ or equivalently, we have
	\begin{equation*}
		L_{\mathcal G}^{\dag} ~ \leq ~ L_{\mathcal P}^{\dag}.
	\end{equation*}
From the linearity property of the trace operator and the fact that $L_{\mathcal P}^{\dag} - L_{\mathcal G}^{\dag}$ is a positive semi-definite matrix, we get
	\begin{eqnarray*}
		\frac{1}{2}\Tr(L_{\mathcal P}^{\dag}-L_{\mathcal G}^{\dag})&=& \frac{1}{2}\Tr(L_{\mathcal P}^{\dag})-	\frac{1}{2}\Tr(L_{\mathcal G}^{\dag}) \nonumber \\ 
		&=&\RR (L_{\mathcal P};M_n) - \RR (L_{\G};M_n) \nonumber \\
		&\geq& 0.
	\end{eqnarray*}
This completes the proof.

\section*{Proof of Theorem \ref{last-thm}}
\label{proof-last-thm}

From our assumptions, we have $L_{\mathcal G_2}  \leq  L_{\mathcal G_1}$, and from the definition, one can verify that 
	\begin{equation*}
		(L_{\mathcal G_1}^{\dag})^2 ~ \leq ~ (L_{\mathcal G_2}^{\dag})^2.
	\end{equation*}
Therefore, using the fact that the trace of a positive semi-definite matrix is always nonnegative, we get
	\begin{eqnarray*}
		\Tr \left ((L_{\mathcal G_2}^{\dag})^2 - (L_{\mathcal G_1}^{\dag})^2 \right) &=& \Tr\left ((L_{\mathcal G_2}^{\dag})^2\right )-\Tr \left ((L_{\mathcal G_1}^{\dag})^2\right ) \\ 
& \geq & 0.
	\end{eqnarray*}
From linearity property of the trace operator, one can conclude that inequality \eqref{subgraph-ineq} holds. The proof of inequality \eqref{subgraph-ineq2} is a direct consequence of Theorem \ref{order-thm}.
\section*{Proof of Corollary \ref{1-coro}}
\label{proof-1-coro}
	By substituting \eqref{propo} in \eqref{sparsity-1}, we get \eqref{tradeoff-1}; and the proof of the additive form inequality is a direct consequence of Theorem \ref{diam-thm} and \eqref{sparsity-1}.

\section*{Proof of Corollary \ref{2-coro}}
\label{proof-2-coro}
	The proof is a direct consequence of Theorems \ref{trace-thm} and \ref{diam-thm}.  

\section*{Proof of Corollary \ref{3-coro}}
\label{proof-3-coro}
	The proof is a direct consequence of Theorem \ref{degree-2-thm} and the definition of $\|\adj\|_{\mathcal{S}_{0,1}}$. 

\section*{Proof of Theorem \ref{13-thm}}
\label{proof-13-thm}

According to \cite{Rojo2000} the following relation holds 
	\begin{equation}
		\lambda_n ~\leq~ \sigma(\G),
		\label{eq1793}
	\end{equation}
where $\sigma(\G)$  is given by (\ref{sss}). By combining inequality (\ref{eq1793}) and (\ref{eq-h}), we get the desired lower bound.

\section*{Proof of Theorem \ref{thm-power-loss}}
\label{proof-power-loss}

From Theorem \ref{2-thm}, we have
	\begin{equation}
		\RR \left ( A_{\G_b}^{(1)};L_{\G_{g}}\oplus \mathbf 0 \right)~=~\frac{1}{2\beta}\mathbf{Tr}(L_{\G_b}^{\dag}L_{\G_{g}}),
		\label{eq_s0}
	\end{equation}
where $L_{\G_b}^{\dag}$ is the Moore-Penrose generalized inverse of the Laplacian matrix $L_{\G_b}$.
According to reference \cite{Gutman2004}, we have 
%
	\begin{equation*}
		L_{\G_b}^{\dag}=-\frac{1}{2}\left (R_{\G_b}-\frac{1}{n}(R_{\G_b}J_{n}+J_{n}R_{\G_b})+\frac{1}{n^2}J_{n}R_{\G_b}J_{n} \right )
	\end{equation*}
where $R_{\G_b}$ is the resistance matrix of the Laplacian matrix $L_{\G_b}$. For a given Laplacian matrix $L_{\G_{g}}$, it is straightforward to verify that $L_{\G_{g}} J_{n} = J_{n} L_{\G_{g}} = 0$. Therefore, we get
	\begin{eqnarray}
		\mathbf{Tr}(L_{\G_b}^{\dag}L_{\G_{g}})&=&
		-\frac{1}{2}\mathbf{Tr}\left ( R_{\G_b} L_{\G_{g}} \right )~=~\sum_{e \in \E_{\mathcal G}}r_e^{(\G_b)}b_e\frac{g_e}{b_e} 				\nonumber \\
		&=& \sum_{e \in \E_{\mathcal G}}\nu_e \alpha_e,
		\label{eq-s}
	\end{eqnarray}
where  $\nu_e=r_e^{(\G_b)}b_e$. From the result of \cite[Lemma 2]{Rapat}, we have that $\sum_{e \in \E_{\mathcal G}}\nu_e =  n-1$. Using this, we can define the weighted mean of the edge parameters $\alpha_e$ for all ${e \in \E_{\mathcal G}}$ as follows
	\begin{eqnarray}
		\bar \alpha=\frac{\sum_{e \in \E_{\mathcal G}}\nu_e \alpha_e}{\sum_{e \in \E_{\mathcal G}}\nu_e}=\frac{\sum_{e \in \E_{\mathcal G}}\nu_e \alpha_e}{n-1}.
		\label{mean}
	\end{eqnarray}
From (\ref{mean}), (\ref{eq-s}) and (\ref{eq_s0}), we conclude the desired result \eqref{total-loss}. 
\section*{Proof of Theorem \ref{edge-trans-thm}}
\label{proof-edge-trans-thm}
 
Similar to the proof of Theorem \ref{thm-power-loss}, we have
	\begin{equation*}
		\RR \left ( A_{\G_b}^{(1)};L_{\G_{g}}\oplus \mathbf 0 \right)~=~\frac{1}{2 \beta}\mathbf{Tr}(L_{\G_b}^{\dag}L_{\G_{g}})~=~\frac{1}{2\beta}\sum_{e \in \E_{\mathcal G}}\nu_e \alpha_e.
	\end{equation*}
Since the coupling graph is  edge-transitive and $\sum_{e \in \E_{\mathcal G}}\nu_e =  n-1$, it follows that $\nu_e=\frac{n-1}{m}$. This completes the proof.

\section*{Proof of Theorem \ref{linear-loss-thm}}
\label{proof-linear-loss-thm}
 
Similar to the proof of Theorem \ref{thm-power-loss}, we have
	\begin{equation*}
		\RR \left ( A_{\G_b}^{(1)};L_{\G_{g}}\oplus \mathbf 0 \right)~=~\frac{1}{2 \beta}\mathbf{Tr}(L_{\mathcal T_b}^{\dag}L_{\mathcal T_{g}})~=~\frac{1}{2\beta}\sum_{e \in \E_{\mathcal T}}\nu_e \alpha_e.
	\end{equation*}
Since the coupling graph is a tree graph, $r_e^{(\G_b)}=b^{-1}_e$ and $\sum_{e \in \E_{\mathcal T}}\nu_e =  n-1$, it follows that $\nu_e=1$.


\newpage

\begin{thebibliography}{99}

\bibitem{Bamieh12}
B.~Bamieh, M.~Jovanovi\'c, P.~Mitra, and S.~Patterson, ``Coherence in
  large-scale networks: Dimension-dependent limitations of local feedback,''
  \emph{IEEE Trans. Autom. Control}, vol.~57, no.~9, pp. 2235--2249, sept.
  2012.

\bibitem{Barooah}
H.~Hao and P.~Barooah, ``Stability and robustness of large platoons of vehicles
  with double-integrator models and nearest neighbor interaction,'' \emph{Int.
  J. Robust and Nonlinear Control}, vol.~23, pp. 2097--2122, 2012.

\bibitem{Siami13siam}
M.~Siami and N.~Motee, ``Robustness and performance analysis of cyclic
  interconnected dynamical networks,'' in \emph{Proc. SIAM Conf. Control and
  Its Appl.}, Jan. 2013, pp. 137--143.

\bibitem{Siami13acc}
M.~Siami, N.~Motee, and G.~Buzi, ``Characterization of hard limits on
  performance of autocatalytic pathways,'' in \emph{Proc. American Control
  Conf.}, 2013, pp. 2313--2318.

\bibitem{Siami13cdc}
M.~Siami and N.~Motee, ``Fundamental limits on robustness measures in networks
  of interconnected systems,'' in \emph{Proc. 52nd IEEE Conf. Decision and
  Control}, Dec. 2013, pp. 67--72.

\bibitem{MoteeCBKD10}
N.~Motee, F.~Chandra, B.~Bamieh, M.~Khammash, and J.~Doyle, ``Performance
  limitations in autocatalytic networks in biology,'' in \emph{Proc. 49th IEEE
  Conf. Decision and Control}, Dec 2010, pp. 4715--4720.

\bibitem{abbas}
W.~Abbas and M.~Egerstedt, ``Robust graph topologies for networked systems,''
  in \emph{Proc. 3rd IFAC Workshop Distributed Estimation and Control in
  Networked Systems}, September 2012, pp. 85--90.

\bibitem{Barooah-acc-12}
H.~Hao and P.~Barooah, ``Improving convergence rate of distributed consensus
  through asymmetric weights,'' in \emph{Proc. American Control Conf.}, 2012,
  pp. 787--792.

\bibitem{Scardovi2010}
L.~Scardovi, M.~Arcak, and E.~Sontag, ``Synchronization of interconnected
  systems with applications to biochemical networks: An input-output
  approach,'' \emph{IEEE Trans. Autom. Control}, vol.~55, no.~6, pp.
  1367--1379, June 2010.

\bibitem{Zelazo2013}
D.~Zelazo, S.~Schuler, and F.~Allg{\"o}wer, ``Performance and design of cycles
  in consensus networks,'' \emph{Syst. Control Lett.}, vol.~62, no.~1, pp.
  85--96, 2013.

\bibitem{Young10}
G.~F. Young, L.~Scardovi, and N.~E. Leonard, ``Robustness of noisy consensus
  dynamics with directed communication,'' in \emph{Proc. American Control
  Conf.}, July 2010, pp. 6312--6317.

\bibitem{Bamieh11}
S.~Patterson and B.~Bamieh, ``Network coherence in fractal graphs,'' in
  \emph{Proc. 50th IEEE Conf. Decision and Control and European Control Conf.},
  Dec. 2011, pp. 6445--6450.

\bibitem{Zelazo-Mesbahi}
D.~Zelazo and M.~Mesbahi, ``Edge agreement: Graph-theoretic performance bounds
  and passivity analysis,'' \emph{IEEE Trans. Autom. Control}, vol.~56, no.~3,
  pp. 544--555, March 2011.

\bibitem{LovisariGarinZampieriResistance}
E.~Lovisari, F.~Garin, and S.~Zampieri, ``Resistance-based performance analysis
  of the consensus algorithm over geometric graphs,'' \emph{SIAM Journal on
  Control and Optimization}, vol.~51, no.~5, pp. 3918--3945, 2013.

\bibitem{Nicola}
N.~Elia, J.~Wang, and X.~Ma, ``Mean square limitations of spatially invariant
  networked systems,'' in \emph{Control of Cyber-Physical Systems}, ser.
  Lecture Notes in Control and Information Sciences, D.~C. Tarraf, Ed.\hskip
  1em plus 0.5em minus 0.4em\relax Springer International Publishing, 2013,
  vol. 449, pp. 357--378.

\bibitem{Lin2014}
F.~Lin, M.~Fardad, and M.~Jovanovic, ``Algorithms for leader selection in
  stochastically forced consensus networks,'' \emph{IEEE Trans. Autom.
  Control}, vol.~59, no.~7, pp. 1789--1802, July 2014.

\bibitem{Jadbabaie13}
A.~Jadbabaie and A.~Olshevsky., ``Combinatorial bounds and scaling laws for
  noise amplification in networks,'' in \emph{Proc. European Control Conf.},
  July 2013, pp. 596--601.

\bibitem{Spielman}
D.~A. Spielman and N.~Srivastava, ``Graph sparsification by effective
  resistances,'' \emph{CoRR}, vol. abs/0803.0929, 2008.

\bibitem{Barooah2006}
P.~Barooah and J.~Hespanha, ``Graph effective resistance and distributed
  control: Spectral properties and applications,'' in \emph{Proc. 45th IEEE
  Conf. Decision and Control}, Dec. 2006, pp. 3479--3485.

\bibitem{Knuth76}
D.~Knuth, ``Big omicron and big omega and big theta,'' \emph{SIGACT News}, pp.
  18--24, 1976.

\bibitem{Bondy}
J.~A. Bondy, \emph{Graph Theory With Applications}.\hskip 1em plus 0.5em minus
  0.4em\relax Elsevier Science Ltd, 1976.

\bibitem{Yu2010}
W.~Yu, G.~Chen, and M.~Cao, ``Some necessary and sufficient conditions for
  second-order consensus in multi-agent dynamical systems,'' \emph{Automatica},
  vol.~46, no.~6, pp. 1089 -- 1095, 2010.

\bibitem{networks-book}
D.~Easley and J.~Kleinberg, \emph{{Networks, Crowds, and Markets: Reasoning
  About a Highly Connected World}}.\hskip 1em plus 0.5em minus 0.4em\relax
  Cambridge, UK: Cambridge University Press, 2010.

\bibitem{Kraning2014}
M.~Kraning, E.~Chu, J.~Lavaei, and S.~Boyd, ``Dynamic network energy management
  via proximal message passing,'' \emph{Foundations and Trends in
  Optimization}, vol.~1, no.~2, pp. 73--126, 2014.

\bibitem{motee-sun-IEEE-TAC-submitted-2014}
N.~Motee and Q.~Sun, ``Sparsity measures for spatially decaying systems,''
  {I}EEE Trans. Autom. Control, submitted.

\bibitem{bamgay13acc}
B.~Bamieh and D.~Gayme, ``The price of synchrony: Resistive losses due to phase
  synchronization in power networks,'' in \emph{Proc. American Control Conf.},
  September 2013, pp. 5815--5820.

\bibitem{5530690}
F.~D\"orfler and F.~Bullo, ``Synchronization and transient stability in power
  networks and non-uniform kuramoto oscillators,'' in \emph{Proc. American
  Control Conf.}, June 2010, pp. 930--937.

\bibitem{Norman}
B.~Norman, \emph{Algebraic Graph Theory}.\hskip 1em plus 0.5em minus
  0.4em\relax Cambridge University Press, Cambridge, 1973.

\bibitem{Mei2004}
Y.~Mei, Y.-H. Lu, Y.~Hu, and C.~Lee, ``Energy-efficient motion planning for
  mobile robots,'' in \emph{IEEE Int. Conf. Robotics and Automation}, vol.~5,
  April 2004, pp. 4344--4349.

\bibitem{marshall11}
A.~W. Marshall, I.~Olkin, and B.~C. Arnold, \emph{Inequalities: Theory of
  Majorization and Its Applications}.\hskip 1em plus 0.5em minus 0.4em\relax
  Springer Science+Business Media, LLC, 2011.

\bibitem{olfati}
R.~Olfati-Saber and R.~Murray, ``Consensus problems in networks of agents with
  switching topology and time-delays,'' \emph{IEEE Trans. Autom. Control},
  vol.~49, no.~9, pp. 1520--1533, Sept 2004.

\bibitem{Rugh}
W.~J. Rugh, \emph{Linear System Theory}.\hskip 1em plus 0.5em minus 0.4em\relax
  Upper Saddle River, NJ, USA: Prentice-Hall, Inc., 1996.

\bibitem{Doyle89}
J.~Doyle, K.~Glover, P.~Khargonekar, and B.~Francis, ``State-space solutions to
  standard $\mathcal{H}_2$ and $\mathcal{H}_{\infty}$ control problems,''
  \emph{IEEE Trans. Autom. Control}, vol.~34, no.~8, pp. 831--847, Aug 1989.

\bibitem{Gutman2003}
I.~Gutman and L.~Pavlovi,``On the coefficients of the laplacian characteristic polynomial of trees,''
  Bulletin. Classe des Sciences
  Math{\'e}matiques et Naturelles. Sciences Math{\'e}matiques, vol. 127,
  no.~28, pp. 31--40, 2003.

\bibitem{tree}
A.~Dobrynin, R.~Entringer, and I.~Gutman, ``Wiener index of trees: Theory and
  applications,'' \emph{Acta Applicandae Mathematica}, vol.~66, no.~3, pp.
  211--249, 2001.

\bibitem{unicyclic}
Y.~Yang and X.~Jiang, ``Unicyclic graphs with extremal kirchhoff index,''
  \emph{MATCH Commun. Math. Comput. Chem}, vol.~60, no.~1, pp. 107--120, 2008.

\bibitem{Yang}
Y.~Yang, ``Bounds for the kirchhoff index of bipartite graphs,'' \emph{J. Appl.
  Math.}, vol. 2012, 2012.

\bibitem{Rapat}
R.~B. Rapat, ``Resistance matrix of a weighted graph,'' \emph{MATCH Commun.
  Math. Comput. Chem.}, vol.~50, pp. 73--82, 2004.

\bibitem{Deng}
H.~Deng, ``On the minimum kirchhoff index of graphs with given number of
  cut-edges,'' \emph{MATCH Commun. Math. Comput. Chem.}, vol.~63, no.~1, pp.
  171--180, 2009.

\bibitem{Zhou2}
B.~Zhou, ``On sum of powers of laplacian eigenvalues and laplacian estrada
  index of graphs,'' 2011, arXiv/1102.1144.

\bibitem{Rojo2000}
O.~Rojo, R.~Soto, and H.~Rojo, ``An always nontrivial upper bound for laplacian
  graph eigenvalues,'' \emph{Linear Algebra and its Appl.}, vol. 312, no. 1--3,
  pp. 155--159, 2000.

\bibitem{Gutman2004}
I.~Gutman and W.~Xiao, ``Generalized inverse of the laplacian matrix and some
  applications,'' \emph{Bulletin. Classe des Sciences Math{\'e}matiques et
  Naturelles. Sciences Math{\'e}matiques}, vol. 129, no.~29, pp. 15--23, 2004.





\end{thebibliography}
\end{document}